\documentclass[12pt]{amsart}

\usepackage{amsmath, amsfonts, amssymb, hyperref}
\usepackage{graphics}
\usepackage[usenames]{color}
\usepackage[all]{xy}
\usepackage{graphicx}
\newcommand{\Q}{\mathcal{Q}} 
 
{\theoremstyle{plain}
\newtheorem{theorem}{Theorem}[section]
\newtheorem{conjecture}[theorem]{Conjecture}
\newtheorem{proposition}[theorem]{Proposition}
\newtheorem{lemma}[theorem]{Lemma}
\newtheorem{corollary}[theorem]{Corollary}
}
{\theoremstyle{definition}
\newtheorem{definition}[theorem]{Definition}
\newtheorem*{definition*}{Definition}
\newtheorem*{proposition*}{Proposition}
\newtheorem*{theorem*}{Result}
\newtheorem{remark}[theorem]{Remark}
\newtheorem{example}[theorem]{Example}

}
\newenvironment{E}{\begin{equation}}{\end{equation}}
\numberwithin{equation}{section}

\newcommand{\coarser}{\succcurlyeq}
\newcommand{\lex}{>_{lex}}
\newcommand{\cell}{\ }
\newcommand{\Qs}{{\mathfrak s}}
\newcommand{\Qr}{{\mathfrak r}}
\newcommand{\tcb}[1]{\textcolor{blue}{#1}}

\newcommand{\tcm}[1]{\textcolor{magenta}{#1}}
\newcommand{\bast}{\tcb{\ast} }
\newcommand{\mast}{\tcm{\ast} }

\newcommand{\ms}{\begin{math}}
\newcommand{\me}{\end{math}}

\newlength\cellsize 
\savebox2{%
\begin{picture}(15,15)
\put(0,0){\line(1,0){15}}
\put(0,0){\line(0,1){15}}
\put(15,0){\line(0,1){15}}
\put(0,15){\line(1,0){15}}
\end{picture}}
\newcommand\cellify[1]{\def\thearg{#1}\def\nothing{}%
\ifx\thearg\nothing
\vrule width0pt height\cellsize depth0pt\else
\hbox to 0pt{\usebox2\hss}\fi%
\vbox to 15\unitlength{
\vss
\hbox to 15\unitlength{\hss$#1$\hss}
\vss}}
\newcommand\tableau[1]{\vtop{\let\\=\cr
\setlength\baselineskip{-16000pt}
\setlength\lineskiplimit{16000pt}
\setlength\lineskip{0pt}
\halign{&\cellify{##}\cr#1\crcr}}}
\savebox3{%
\begin{picture}(15,15)
\put(0,0){\line(1,0){15}}
\put(0,0){\line(0,1){15}}
\put(15,0){\line(0,1){15}}
\put(0,15){\line(1,0){15}}
\end{picture}}
\newcommand\expath[1]{%
\hbox to 0pt{\usebox3\hss}%
\vbox to 15\unitlength{
\vss
\hbox to 15\unitlength{\hss$#1$\hss}
\vss}}

\newcommand{\comment}[1]{\vspace{5 mm}\par \noindent
\marginpar{\textsc{Comment}}
\framebox{\begin{minipage}[c]{0.95 \textwidth}
 #1 \end{minipage}}\vspace{5 mm}\\}

\renewcommand{\comment}[1]{}

\begin{document}
\title[Equality of Schur Q-functions]{Composition of transpositions and equality of ribbon Schur $Q$-functions}

\author{Farzin Barekat}
\address{Department of Mathematics, University of British Columbia, Vancouver, BC V6T 1Z2, Canada}
\email{}

\author{Stephanie van Willigenburg}
\address{Department of Mathematics, University of British Columbia, Vancouver, BC V6T 1Z2, Canada}
\email{steph@math.ubc.ca}


\maketitle

\begin{abstract}
We introduce a new operation on skew diagrams called composition of transpositions, and use it and a Jacobi-Trudi style formula to derive equalities on skew Schur $Q$-functions whose indexing shifted skew diagram is an ordinary skew diagram. When this skew diagram is a ribbon, we conjecture necessary and sufficient conditions for equality of ribbon Schur $Q$-functions. Moreover, we determine all relations between ribbon Schur $Q$-functions; show they supply a $\mathbb{Z}$-basis for skew Schur $Q$-functions; assert their irreducibility; and show that the non-commutative analogue of ribbon Schur $Q$-functions is the flag $h$-vector of Eulerian posets.
\end{abstract}
\tableofcontents

\section{Introduction}\label{sec:intro} 
In the algebra of symmetric functions there is interest in determining when two skew Schur functions are equal \cite{HDL, gut, HDL3, HDL2, vW}. The equalities are described in terms of equivalence relations on skew diagrams. It is consequently natural to investigate whether new equivalence relations on skew diagrams arise when we restrict our attention to the subalgebra of skew Schur $Q$-functions. This is a particularly interesting subalgebra to study since the combinatorics of skew Schur $Q$-functions also arises in the representation theory of the twisted symmetric group \cite{Bess, ShawvW, StemP}, and the theory of enriched $P$-partitions \cite{Stem}, and hence skew Schur $Q$-function equality would impact these areas. The study of skew Schur $Q$-function equality was begun in \cite{QEq}, where a series of technical conditions classified when a skew Schur $Q$-function is equal to a Schur $Q$-function. In this paper we extend this study to the equality of ribbon Schur $Q$-functions. Our motivation for focussing on this family is because the study of ribbon Schur function equality is fundamental to the general study of skew Schur function equality, as evidenced by \cite{HDL, HDL3, HDL2}. Our method of proof is to study a slightly more general family of skew Schur $Q$-functions, and then restrict our attention to ribbon Schur $Q$-functions. Since the combinatorics of skew Schur $Q$-functions is more technical than that of skew Schur functions, we provide detailed proofs to highlight the subtleties needed to be considered for the general study of equality of skew Schur $Q$-functions.
The rest of this paper is structured as follows.

In the next section we review operations on skew diagrams, introduce the skew diagram operation \emph{composition of transpositions} and derive some basic properties for it, including associativity in Proposition~\ref{prop:assoc}. In Section~\ref{sec:schurq}  we recall $\Omega$, the algebra of Schur $Q$-functions, discover new bases for this algebra in Proposition~\ref{prop:sbasis} and Corollary~\ref{cor:rbasis}. We see the prominence of ribbon Schur $Q$-functions in the latter, which states

\begin{theorem*} The set of all ribbon Schur $Q$-functions $\Qr _\lambda$, indexed by strict partitions $\lambda$, forms a $\mathbb{Z}$-basis for $\Omega$.
\end{theorem*}

Furthermore we determine all relations between ribbon Schur $Q$-functions in Theorems~\ref{ribbonrelations} and \ref{ribbonrelations2}. The latter is particularly succinct:

\begin{theorem*}All relations amongst ribbon Schur $Q$-functions are generated by the multiplication rule $\Qr _\alpha \Qr _\beta = \Qr _{\alpha \cdot \beta} + \Qr _{\alpha \odot \beta}$ for compositions $\alpha, \beta$, and $\Qr _{2m} = \Qr _{1^{2m}}$ for $m\geq 1$.
\end{theorem*}

In Section~\ref{sec:eqskewschurq} we determine a number of instances when two ordinary skew Schur $Q$-functions are equal including a necessary and sufficient condition in Proposition~\ref{prop:power2}. Our main theorem on equality is Theorem~\ref{the:bigone}, which  is dependent on composition of transpositions denoted $\bullet$, transposition denoted $^t$, and antipodal rotation denoted $^\circ$:

\begin{theorem*}
For ribbons $\alpha _1, \ldots , \alpha _m$ and skew diagram $D$ the ordinary skew Schur $Q$-function indexed by
$$\alpha _1 \bullet \cdots \bullet \alpha _m \bullet D$$is equal to the ordinary skew Schur $Q$-function indexed by
$$\beta _1 \bullet \cdots \bullet \beta _m \bullet E$$where
$$ 
\beta _i \in \{ \alpha _i, \alpha _i ^t, \alpha _i ^\circ , (\alpha _i ^t)^\circ = (\alpha _i ^\circ)^t\} \quad 1\leq i \leq m,$$
$$ 
E\in \{ D, D^t, D^\circ , (D^t)^\circ = (D^\circ)^t\}.$$
\end{theorem*}
 We restrict our attention  to ribbon Schur $Q$-functions again in Section~\ref{sec:ribbons}, and derive further ribbon specific properties  including irreducibility in Proposition~\ref{prop:irrrib}, and that the non-commutative analogue of ribbon Schur $Q$-functions is the flag $h$-vector of Eulerian posets in Theorem~\ref{the:commconnection}.

\section*{Acknowledgements}\label{sec:ack}
The authors would like to thank Christine Bessenrodt, Louis Billera and Hugh Thomas for helpful conversations,  Andrew Rechnitzer for programming assistance, and the referee for helpful comments. John Stembridge's QS package helped to generate the pertinent data. Both authors were supported in part by the National Sciences and Engineering Research Council of Canada.

\section{Diagrams}\label{sec:diagrams}
A \emph{partition}, $\lambda$, of a positive integer $n$, is a list of positive integers $\lambda _1  \geq \cdots \geq \lambda _k >0$ whose sum is $n$. We denote this by $\lambda \vdash n$, and for convenience denote the empty partition of 0 by 0. We say that a partition is \emph{strict} if $\lambda _1 > \cdots > \lambda _k >0$. If we remove the weakly decreasing criterion from the partition definition, then we say the list is a composition. That is, a \emph{composition}, $\alpha$, of a positve integer $n$ is a list of positive integers $\alpha _1 \cdots \alpha _k$ whose sum is $n$. We denote this by $\alpha \vDash n$. Notice that any composition $\alpha = \alpha _1 \cdots \alpha _k$ determines a partition, denoted $\lambda (\alpha)$, where $\lambda (\alpha)$ is obtained by reordering $\alpha _1 , \ldots , \alpha _k$ in weakly decreasing order. Given a composition $\alpha = \alpha _1 \cdots \alpha _k\vDash n$ we call the $\alpha _i$ the \emph{parts}  of $\alpha$, $n= :|\alpha |$ the \emph{size} of $\alpha$ and $k=: \ell (\alpha)$ the \emph{length} of $\alpha$. There also exists three partial orders on compositions, which will be useful to us later.

Firstly, given two compositions $\alpha = \alpha _1\cdots \alpha _{\ell (\alpha)}$, $\beta = \beta _1 \cdots \beta _{\ell(\beta)} \vDash n$ we say $\alpha $ is a \emph{coarsening} of $\beta$ (or $\beta$ is a \emph{refinement} of $\alpha$), denoted $\alpha \coarser \beta$ if adjacent parts of $\beta$ can be added together to yield the parts of $\alpha$, for example, $5312 \coarser 1223111$. Secondly, we say $\alpha$ \emph{dominates} $\beta$, denoted $\alpha \geq \beta$ if
$\alpha _1 + \cdots +\alpha _i \geq \beta _1 + \cdots + \beta _i $ for $i=1, \ldots , \min\{\ell (\alpha), \ell(\beta)\}.$ Thirdly, we say $\alpha$ is \emph{lexicographically greater than} $\beta$, denoted $\alpha \lex \beta$ if $\alpha \neq \beta$ and the first $i$ for which $\alpha _i \neq \beta _i$ satisfies $\alpha _i > \beta _i$.

From partitions we can also create diagrams as follows. Let $\lambda $ be a partition. Then the array of left justified cells containing $\lambda _i$ cells in the $i$-th row from the top is called the \emph{(Ferrers or Young) diagram} of $\lambda$, and we abuse notation by also denoting it by $\lambda$.
Given two diagrams $\lambda, \mu$ we say $\mu$ is \emph{contained} in $\lambda$, denoted $\mu \subseteq \lambda$ if $\mu _i \leq \lambda _i$ for all $i=1, \ldots ,\ell(\mu)$. Moreover, if $\mu \subseteq \lambda$ then the \emph{skew diagram} $D=\lambda /\mu$ is obtained from the diagram of $\lambda$ by removing the diagram of $\mu$ from the top left corner.
The \emph{disjoint union} of two skew diagrams $D_1$ and $D_2$, denoted $D_1 \oplus D_2$, is obtained by placing $D_1$ strictly north and east of $D_2$ such that $D_1$ and $D_2$ occupy no common row or column.
We say a skew diagram is \emph{connected} if it cannot be written as $D_1\oplus D_2$ for two non-empty skew diagrams $D_1, D_2$. If a connected skew diagram additionally contains no $2\times 2$ subdiagram then we call it a \emph{ribbon}. Ribbons will be an object of focus for us later, and hence for ease of referral we now recall the well-known correspondence between ribbons and compositions. Given a ribbon with $\alpha _1$ cells in the 1st row, $\alpha _2$ cells in the 2nd row, $\ldots$, $\alpha _{\ell(\alpha)}$ cells in the last row, we say it corresponds to the composition $\alpha _1 \cdots \alpha _{\ell(\alpha)}$, and we abuse notation by denoting the ribbon by $\alpha$ and noting it has $|\alpha|$ cells.

\begin{example}
$\lambda / \mu = 3221 / 11  = \tableau{&\cell&\cell\\&\cell\\ 
\cell&\cell \\\cell} = 2121 = \alpha$. \end{example}

\subsection{Operations on diagrams}\label{subsec:ops}
In this subsection we introduce operations on skew diagrams that will enable us to describe more easily when two skew Schur $Q$-functions are equal. We begin by recalling three classical operations: transpose, antipodal rotation, and shifting.

Given a diagram $\lambda = \lambda _1 \cdots \lambda _{\ell(\lambda)}$ we define the \emph{transpose} (or \emph{conjugate}), denoted $\lambda ^t$, to be the diagram containing $\lambda _i$ cells in the $i$-th column from the left. We extend this definition to skew diagrams by defining the transpose of $\lambda /\mu$ to be $(\lambda /\mu)^t:=\lambda ^t/\mu ^t$ for diagrams $\lambda, \mu$. Meanwhile, the \emph{antipodal rotation} of $\lambda /\mu$, denoted $(\lambda /\mu )^\circ$, is obtained by rotating $\lambda /\mu$ 180 degrees in the plane. Lastly, if $\lambda , \mu$ are strict partitions then we define the \emph{shifted} skew diagram of $\lambda /\mu$, denoted $(\widetilde{\lambda/\mu})$, to be the array of cells obtained from $\lambda /\mu$ by shifting the $i$-th row from the top $(i-1)$ cells to the right for $i>1$.

\begin{example} If $\lambda = 5421, \mu = 31$ then
$$\lambda/\mu = \tableau{&&&\cell&\cell\\
&\cell&\cell&\cell\\
\cell&\cell\\
\cell}, (\lambda /\mu)^t = \tableau{&&\cell &\cell \\
&\cell&\cell\\
& \cell\\
\cell&\cell\\
\cell}, (\lambda /\mu )^\circ = \tableau{&&&&\cell\\
&&&\cell &\cell\\
& \cell&\cell&\cell\\
\cell&\cell}, (\widetilde{\lambda/\mu}) = \tableau{&\cell&\cell\\
\cell&\cell&\cell\\
\cell&\cell\\
&\cell}.$$
\end{example}

We now recall three operations that are valuable in describing when two skew Schur functions are equal, before introducing a new operation. The first two operations, concatenation and near concatenation, are easily obtained from the disjoint union of two skew diagrams $D_1, D_2$. Given $D_1\oplus D_2$ their \emph{concatenation} $D_1\cdot D_2$ (respectively, \emph{near concatenation} $D_1 \odot D_2$) is formed by moving all the cells of $D_1$ exactly one cell west (respectively, south).

\begin{example}If $D_1 = 21, D_2=32$ then
$$D_1\oplus D_2 = \tableau{&&&\cell &\cell\\
&&& \cell\\
\cell&\cell&\cell\\
\cell&\cell}\ ,\quad D_1 \cdot D_2 = \tableau{&&\cell&\cell\\
&&\cell\\
\cell&\cell&\cell\\
\cell&\cell}\ , \quad D_1\odot D_1 = \tableau{&&&\cell&\cell\\
\cell&\cell&\cell&\cell\\
\cell&\cell}\ .$$
\end{example}

For the third operation recall that $\cdot$ and $\odot$ are each associative and associate with each other \cite[Section 2.2]{HDL2} and hence any string of operations on diagrams $D_1, \ldots , D_k$
$$D_1\bigstar _1 D_2 \bigstar _2 \cdots \bigstar _{k-1} D_k$$
in which each $\bigstar _i$ is either $\cdot$ or $\odot$ is well-defined without parenthesization. Also recall from \cite{HDL2} that a ribbon with $|\alpha |=k$ can be uniquely written as
$$\alpha =\square\bigstar _1 \square \bigstar _2 \cdots \bigstar _{k-1} \square$$where $\square$ is the diagram with one cell. Consequently, given a composition $\alpha$ and skew diagram $D$ the operation \emph{composition of compositions} is
$$\alpha \circ D = D\bigstar _1 D\bigstar _2 \cdots \bigstar _{k-1} D.$$This third operation was introduced in this way in \cite{HDL2} and we modify this description to define our fourth, and final, operation \emph{composition of transpositions} as
\begin{equation} \alpha \bullet D=\left \{ \begin{array}{ll} 
                 D\bigstar _1 D^t\bigstar _2 D \bigstar _3 D ^t \cdots \bigstar _{k-1} D & \hbox{if } |\alpha| \hbox{ is odd } \\
                 D\bigstar _1 D^t\bigstar _2 D \bigstar _3 D ^t \cdots \bigstar _{k-1} D^t & \hbox{if } |\alpha| \hbox{ is even. } 
                \end{array} \right.  \label{eq:compoftrans} \end{equation} 
                
We refer to $\alpha \circ D$ and $\alpha \bullet D$ as  consisting of blocks of $D$ when we wish to highlight the  dependence on $D$.

\begin{example}
Considering our block to be $D = 31$ and using coloured $\ast$ to highlight the blocks
$$21 \circ D = \tableau{&&&&&\bast&\bast&\bast\\
&&\mast&\mast&\mast&\bast\\
&&\mast\\
\bast&\bast&\bast\\
\bast}\ \mbox{ and }\ 21\bullet D = \tableau{&&&&\bast&\bast&\bast\\
&&\mast&\mast&\bast\\
&&\mast\\
&&\mast\\
\bast&\bast&\bast\\
\bast
}\ .$$Observe that if we consider the block $D=2$, then the latter ribbon can also be described as $312 \bullet 2$:
$$\tableau{&&&&\mast&\bast&\bast\\
&&\bast&\bast&\mast\\
&&\mast\\
&&\mast\\
\mast&\bast&\bast\\
\mast
}\ .$$
\end{example}

This last operation will be the focus of our results, and hence we now establish some of its basic properties.

\subsection{\texorpdfstring{Preliminary properties of $\bullet$}{Preliminary properties of bullet}}\label{subsec:bulletprops}
Given a ribbon $\alpha$ and skew diagram $D$ it is straightforward to verify using \eqref{eq:compoftrans} that
\begin{equation} ( \alpha \bullet D)^\circ=\left \{ \begin{array}{ll} 
                 \alpha^\circ \bullet D^\circ & \hbox{if } |\alpha| \hbox{ is odd } \\
                 \alpha^\circ \bullet (D^t)^\circ & \hbox{if } |\alpha| \hbox{ is even } 
                \end{array} \right.  \label{dirotation} \end{equation} 
and
\begin{equation} ( \alpha \bullet D)^t=\left \{ \begin{array}{ll} 
  \alpha^t \bullet D^t & \hbox{if } |\alpha| \hbox{ is odd } \\
                 \alpha^t \bullet D & \hbox{if } |\alpha| \hbox{ is even. } 
                \end{array} \right.  \label{ditransposition} \end{equation} 
                
We can also verify that $\bullet$ satisfies an associativity property, whose proof illustrates some of the subtleties of $\bullet$.

\begin{proposition}\label{prop:assoc}
Let $\alpha, \beta$ be ribbons and $D$ a skew diagram. Then
$$\alpha \bullet (\beta \bullet D) = (\alpha \bullet \beta) \bullet D.$$
\end{proposition}

\begin{proof}
First notice that, if we decompose the $\beta \bullet  D$ components of $\alpha \bullet (\beta \bullet  D)$ into blocks of $ D$ then the $ D$ blocks are alternating in appearance as $D$ or $D^t$ as is in $(\alpha \bullet \beta) \bullet  D$. Furthermore both $\alpha \bullet (\beta \bullet  D)$ and $(\alpha \bullet \beta) \bullet  D$ are comprised of $|\alpha|\times|\beta|$ blocks of $ D$. The only remaining thing is to show that the $i$-th and $i+1$-th block of $ D$ are joined in the same manner (i.e. near concatenated or concatenated) in both $\alpha \bullet (\beta \bullet  D)$ and $(\alpha \bullet \beta) \bullet  D$.
\\
For a ribbon $\gamma$ let 
$$f^{\gamma}(i)=\left \{ \begin{array}{ll} 
                 -1 & \hbox{if in the ribbon $\gamma$, the $i$-th and $i+1$-th cell are near concatenated}  \\
                 1 & \hbox{if in the ribbon $\gamma$, the $i$-th and $i+1$-th cell are  concatenated.}
                \end{array} \right. $$
{\it Case 1: i=$|\beta|q$.} Note that $\beta \bullet  D$ has $|\beta|$ blocks of $ D$. Therefore, the way that the $i$-th and $i+1$-th blocks of $ D$ are joined in $\alpha \bullet (\beta \bullet  D)$ is given by $f^{\alpha}(q)$. Now in $(\alpha \bullet \beta) \bullet  D$ the way that the $i$-th and $i+1$-th blocks of $ D$ are joined is given by $f^{\alpha \bullet \beta}(i)$, which is equal to $f^{\alpha}(q)$. 
\\
{\it Case 2: i=$|\beta|q+r$ where $r\neq 0$.} Note that $f^{\gamma^t}(i)=-f^{\gamma}(i)$. Since in $\alpha \bullet \beta$, the $\beta$ components are alternating in appearance as $\beta$, $\beta ^t$, the way that the $i$-th and $i+1$-th block of $ D$ are joined in $(\alpha \bullet \beta) \bullet  D$ is given by $f^{\alpha\bullet\beta}(i)=(-1)^{q}f^{\beta}(r)$. For $\alpha \bullet (\beta \bullet  D)$, note that the $i$-th and $i+1$-th blocks of $ D$ are part of $\beta \bullet  D$, hence they are joined given by $(-1)^{q}f^{\beta}(r)$, where $(-1)^{q}$ comes from the fact that we are using $\beta \bullet  D$ and its transpose alternatively to form $\alpha \bullet (\beta \bullet  D)$.
\end{proof}

\section{\texorpdfstring{Skew Schur $Q$-functions}{Skew Schur Q-functions}}\label{sec:schurq} We now introduce our objects of study, skew Schur $Q$-functions. Although they can be described in terms of Hall-Littlewood functions at $t=-1$ we define them combinatorially for later use.

Consider the alphabet
$$1'<1<2'<2<3'<3 \cdots .$$Given a shifted skew diagram $(\widetilde{\lambda /\mu})$ we define a \emph{weakly amenable tableau}, $T$, of \emph{shape} $(\widetilde{\lambda /\mu})$ to be a filling of the cells of $(\widetilde{\lambda /\mu})$ such that
\begin{enumerate}
\item the entries in each row of $T$ weakly increase
\item the entries in each column of $T$ weakly increase
\item each row contains at most one $i'$ for each $i\geq 1$
\item each column contains at most one $i$ for each $i\geq 1$.
\end{enumerate}
We define the {content} of $T$ to be
$$c(T)=c_1(T)c_2(T)\cdots$$where
$$c_i(T)= |\ i\ | +|\ i'\ |$$
and $|\ i\ |$ is the number of times $i$ appears in $T$, whilst $|\ i'\ |$ is the number of times $i'$ appears in $T$. The monomial associated to $T$ is given by
$$x^T:=x_1 ^{c_1(T)}x_2 ^{c_2(T)}\cdots $$and the \emph{skew Schur $Q$-function}, $Q_{\lambda /\mu}$, is then
$$Q_{\lambda /\mu} = \sum _T x^T$$where the sum is over all weakly amenable tableau $T$ of shape $(\widetilde{\lambda /\mu})$.
Two skew Schur $Q$-functions that we will be particularly interested in are ordinary skew Schur $Q$-functions and ribbon Schur $Q$-functions.

If $(\widetilde{\lambda /\mu}) = D$ where $D$ is a skew diagram then we define
$$\Qs _D:= Q_{\lambda /\mu}$$and call it an \emph{ordinary skew Schur $Q$-function}. If, furthermore, $(\widetilde{\lambda /\mu})$ is a ribbon, $\alpha$, then we define
$$\Qr _\alpha := Q_{\lambda /\mu}$$and call it a \emph{ribbon Schur $Q$-function}. 

Skew Schur $Q$-functions lie in the algebra $\Omega$, where
$$\Omega = \mathbb{Z} [ q_1, q_2, q_3, \ldots ] \equiv \mathbb{Z} [ q_1, q_3, q_5, \ldots ]$$and
$q_n = Q_n$. The $q_n$ satisfy
\begin{equation}\sum _{r+s = n} (-1)^r q_rq_s = 0, \label{eq:qrels}\end{equation}which will be useful later, but for now note that for any set of countable indeterminates $x_1, x_2, \ldots$ the expression $\sum _{r+s = n} (-1)^r x_rx_s$ is often denoted $\chi _n$ and is called the \emph{$n$-th Euler form}.

Moreover, if $\lambda = \lambda _1 \cdots \lambda _{\ell(\lambda)}$ is a partition and we define
$$q_\lambda := q_{\lambda _1}\cdots q_{\ell(\lambda)},\quad q_0=1$$then

\begin{proposition}\cite[8.6(ii)]{MacD}\label{prop:qbasis}
The set $\{ q_\lambda \} _{\lambda \vdash n\geq 0}$, for $\lambda$ strict, forms a $\mathbb{Z}$-basis of $\Omega$.
\end{proposition}

This is not the only basis of $\Omega$ as we will see in Proposition~\ref{prop:sbasis}.

\subsection{\texorpdfstring{Symmetric functions and $\theta$}{Symmetric functions and theta}}\label{subsec:symmap}
It transpires that the $\Qs _D$ and $\Qr _\alpha$ can also be obtained from symmetric functions.
Let $\Lambda$ be the subalgebra of $\mathbb{Z}[x_1, x_2,  \ldots]$ with countably many variables $x_1, x_2,  \ldots$
 given by $\Lambda = \mathbb{Z}[e_1, e_2, \ldots ] = \mathbb{Z}[h_1, h_2, \ldots ]$ where
$e_n = \prod _{ i_1 < \cdots <i_n} x_{i_1}\cdots x_{i_n}$ is the \emph{$n$-th elementary symmetric function} and 
$h_n = \prod _{ i_1 \leq \cdots \leq i_n} x_{i_1}\cdots x_{i_n}$ is the \emph{$n$-th homogeneous symmetric function}. Moreover,  if $\lambda = \lambda _1\cdots \lambda _{\ell(\lambda)}$ is a partition and we define
$e_\lambda := e_{\lambda _1}\cdots e_{\ell(\lambda)}$, 
$h_\lambda := h_{\lambda _1}\cdots h_{\ell(\lambda)}$, and
$e_0=h_0=1$ then

\begin{proposition}\cite[I.2]{MacD}\label{prop:ehbasis}
The sets $\{ e_\lambda \} _{\lambda \vdash n\geq 0}$ and $\{ h_\lambda \} _{\lambda \vdash n\geq 0}$, each form a $\mathbb{Z}$-basis of $\Lambda$.
\end{proposition} 

Given a skew diagram, $\lambda /\mu$ we can use the \emph{Jacobi-Trudi determinant formula} to describe the \emph{skew Schur function} $s_{\lambda /\mu}$ as
\begin{equation}\label{eq:jth}s_{\lambda /\mu} = \det (h _{\lambda _i-\mu _j -i+j}) _{i,j = 1} ^{\ell(\lambda)}\end{equation}and via the involution $\omega:\Lambda \rightarrow \Lambda$ mapping $\omega (e_n)=h_n$ we can deduce
\begin{equation}\label{eq:jte}s_{(\lambda /\mu)^t} = \det (e _{\lambda _i-\mu _j -i+j}) _{i,j = 1} ^{\ell(\lambda)}\end{equation}where
$\mu_i = 0, i>\ell(\mu)$ and $h_n=e_n=0$ for  $n<0$.

If, furthermore, $\lambda/\mu$ is a ribbon $\alpha$ then we define
$$r_\alpha := s _{\lambda /\mu}$$and call it a \emph{ribbon Schur function}.

To obtain an algebraic description of our ordinary and ribbon Schur $Q$-functions we need the graded surjective ring homomorphism 
$$\theta : \Lambda \longrightarrow \Omega$$that  satisfies \cite{Stem}
$$\theta (h_n)=\theta (e_n)=q_n, \quad \theta (s_D) = \Qs _D,\quad \theta(r_\alpha )=\Qr _\alpha$$for any skew diagram $D$ and ribbon $\alpha$. The homomorphism $\theta$ enables us to immediately determine a number of properties of ordinary skew and ribbon Schur $Q$-functions.

\begin{proposition} Let $\lambda /\mu$ be a skew diagram and $\alpha $ a ribbon. Then
\begin{equation}\label{eq:Qrot}
\Qs _{\lambda /\mu} = \Qs _{(\lambda /\mu)^\circ}
\end{equation}
\begin{equation}\label{eq:Qtr}
\Qs _{\lambda /\mu} = \det (q _{\lambda _i-\mu _j -i+j}) _{i,j = 1} ^{\ell(\lambda)} = \Qs _{(\lambda /\mu)^t}
\end{equation}
\begin{equation}\label{eq:Qrib}
\Qr _\alpha = (-1)^{\ell(\alpha)} \sum _{\beta \coarser \alpha} (-1) ^{\ell(\beta)} q _{\lambda (\beta)}.
\end{equation}
Moreover, for $D,E$ being skew diagrams and $\alpha, \beta$ being ribbons
\begin{equation}\label{eq:Qmult}
\Qs _D\Qs _E = \Qs _{D\cdot E} + \Qs _{D\odot E}
\end{equation}
\begin{equation}\label{eq:Qribmult}
\Qr _\alpha\Qr _\beta = \Qr _{\alpha\cdot \beta} + \Qr _{\alpha\odot \beta}.
\end{equation}
\end{proposition}

\begin{proof}
The first equation follows from applying $\theta$ to \cite[Exercise 7.56(a)]{ECII}. The second equation follows from applying $\theta$ to \eqref{eq:jth} and \eqref{eq:jte}. The third equation follows from applying $\theta$ to \cite[Proposition 2.1]{HDL}. The fourth and fifth equations follow from applying $\theta$ to  \cite[Proposition 4.1]{HDL2} and \cite[(2.2)]{HDL}, respectively.
\end{proof}

\subsection{\texorpdfstring{New bases and relations in $\Omega$}{New bases and relations in Omega}}\label{subsec:basesandrels}
The map $\theta$ is also useful for describing bases for $\Omega$ other than the basis given in Proposition~\ref{prop:qbasis}.

\begin{definition}
If $D$ is a skew diagram, then let $srl(D)$ be the partition determined by the (multi)set of row lengths of $D$.
\end{definition}

\begin{example}
$$D = \tableau{ &\cell&\cell\\
\cell&\cell&\cell\\
\cell&\cell\\
\cell}\  \quad srl(D) = 3221$$
\end{example}

\begin{proposition} \label{prop:sbasis}
Let $\mathfrak{D}$ be a set of skew diagrams such that for all $D\in \mathfrak{D}$ we have $srl(D)$ is a strict partition, and for all strict partitions $\lambda$ there exists exactly one $D\in \mathfrak{D}$ satisfying $srl(D)=\lambda$.
Then the set $\{ \Qs _D \} _{D\in \mathfrak{D}}$ forms a $\mathbb{Z}$-basis of $\Omega$.
\end{proposition}

\begin{proof}
Let $D$ be any skew diagram such that $srl(D)=\lambda$. By \cite[Proposition 6.2(ii)]{HDL2}, we know that $h_\lambda$ has the lowest subscript in dominance order when we expand the skew Schur function $s_D$ in terms of complete symmetric functions. That is 
$$s_D=h_\lambda+\hbox{\scriptsize a sum of $h_\mu$'s where $\mu$ is a partition with $\mu>\lambda$}.$$
Now applying $\theta$ to this equation and using \cite[(8.4)]{MacD}, we conclude that
\begin{equation} {\mathfrak s}_D=q_\lambda+\hbox{\scriptsize a sum of $q_\mu$'s where $\mu$ is a strict partition with $\mu>\lambda$} \label{sdtoq}.\end{equation}
Hence by Proposition~\ref{prop:qbasis}, the set of ${\mathfrak s}_D$, $D\in{\mathfrak D}$, forms a basis of $\Omega$.

The equation \eqref{sdtoq} implies that if we order $\lambda$'s and $srl(D)$'s in lexicographic order the transition matrix that takes ${\mathfrak s}_D$'s to $q_\lambda$'s is unitriangular with integer coefficients.  Thus,  the transition matrix that takes  $q_\lambda$'s to ${\mathfrak s}_D$'s  is unitriangular with integer coefficients. Hence 
\begin{equation}q_\lambda={\mathfrak s}_D+ \hbox{\scriptsize a sum of ${\mathfrak s}_E$'s where $srl(E)$ is a strict partition and $srl(E)>srl(D)$} \label{qtosd}\end{equation}
where $E,D \in {\mathfrak D}$ and $srl(D)=\lambda$. 

Combining Proposition~\ref{prop:qbasis} with \eqref{qtosd} it follows  that the set of ${\mathfrak s}_D$, $D\in {\mathfrak D}$, forms a ${\mathbb Z}$-basis of $\Omega$.
\end{proof}

\begin{corollary}\label{cor:rbasis}
The set $\{ \Qr _\lambda \} _{\lambda \vdash n \geq 0}$, for $\lambda$ strict, forms a $\mathbb{Z}$-basis of $\Omega$.
\end{corollary}

We can now describe a set of relations that generate \emph{all} relations amongst ribbon Schur $Q$-functions.

\begin{theorem} 
\label{ribbonrelations} 
Let {$z_{\alpha}, \alpha \vDash n, n\ge 1$} be commuting indeterminates.  Then 
as algebras, ${\Omega}$ is isomorphic to the quotient 
$${\Q[z_{\alpha}]/ 
 \langle z_{\alpha}\ z_{\beta}-z_{\alpha\cdot \beta} - z_{\alpha \odot \beta}, \chi _2, \chi_4, \ldots\rangle}$$where
 $\chi_{2m}$ is the even Euler form $\chi _{2m} = \sum _{r+s = 2m} (-1)^{r} z_rz_s$. Thus, all relations amongst ribbon Schur $Q$-functions are generated by $\Qr_{\alpha}\ \Qr_{\beta}= \Qr_{\alpha\cdot \beta} + \Qr_{\alpha \odot \beta}$ and  $\sum _{r+s = 2m} (-1)^{r} \Qr_r\Qr_s = 0$, $m\geq 1$.
\end{theorem} 

\begin{proof} Consider the map $\varphi:\Q[z_{\alpha}] \rightarrow \Omega$ defined by 
 $z_{\alpha} \mapsto {\mathfrak r}_\alpha$.  This map is surjective since the ${\mathfrak r}_{\alpha}$ 
generate $\Omega$ by Corollary~\ref{cor:rbasis}.   Grading $\Q[z_{\alpha}]$ by setting the degree of $z_{\alpha}$ 
to be $n=|\alpha|$ makes $\varphi$ homogeneous. 
To see that $\varphi$ induces an isomorphism with the quotient, note that 
$\Q[z_{\alpha}]/ 
 \langle z_{\alpha}\ z_{\beta}-z_{\alpha\cdot \beta} - z_{\alpha \odot \beta},\chi _2, \chi_4, \ldots\rangle$ 
 maps onto $\Q[z_{\alpha}]/\ker \varphi \simeq\Omega $, 
 since $ \langle z_{\alpha}\ z_{\beta}-z_{\alpha\cdot \beta} - z_{\alpha \odot \beta}, \chi _2, \chi_4,\ldots\rangle 
 \subset \ker\varphi$ as we will see below. 
  
 It then suffices to show that  the degree $n$ component of 
$$\Q[z_{\alpha}]/ 
 \langle z_{\alpha}\ z_{\beta}-z_{\alpha\cdot \beta} - z_{\alpha \odot \beta}, \chi _2, \chi_4,\ldots \rangle$$is generated by 
the images of the $z_{\lambda}, \lambda \vdash n$, $\lambda$ is a strict partition, and so has dimension at 
most the number of  partitions of $n$ with distinct parts.  

We show $ \langle z_{\alpha}\ z_{\beta}-z_{\alpha\cdot \beta} - z_{\alpha \odot \beta}, \chi _2, \chi_4,\ldots\rangle 
 \subset \ker\varphi$ as follows. 

From \cite[p 251]{MacD} we know that 
\begin{E} 2q_{2x}=q_{2x-1}q_1-q_{2x-2}q_2+\cdots+q_1q_{2x-1} \label{qrelations} \end{E}and
since $q_i={\mathfrak r}_i$, we can rewrite the above equation 
$$2{\mathfrak r}_{2x}={\mathfrak r}_{2x-1}{\mathfrak r}_1-{\mathfrak r}_{2x-2}{\mathfrak r}_2+\cdots+{\mathfrak r}_1{\mathfrak r}_{2x-1}.$$
Substituting ${\mathfrak r}_{2x-i}{\mathfrak r}_i={\mathfrak r}_{2x}+{\mathfrak r}_{(2x-i)i}$ and simplifying, we get
$${\mathfrak r}_{2x}={\mathfrak r}_{(2x-1)1}-{\mathfrak r}_{(2x-2)2}+\cdots+(-1)^{x+1}{\mathfrak r}_{xx}+\cdots+{\mathfrak r}_{1(2x-1)}.$$Together with \eqref{eq:Qribmult} we have $ \langle z_{\alpha}\ z_{\beta}-z_{\alpha\cdot \beta} - z_{\alpha \odot \beta}, \chi _2, \chi_4,\ldots\rangle 
 \subset \ker\varphi$.

Now we show that if we have the following relations then every $z_\gamma$ can be written as the sum of $z_\lambda$'s where the $\lambda$'s are strict partitions. 
\begin{E} 
\left \{ \begin{array}{l} 
                 z_\alpha z_\beta=z_{\alpha\cdot\beta}+z_{\alpha\odot\beta}  \\
                  z_2=z_{11}\\
                  z_4=z_{31}-z_{22}+z_{13}\\
                  \vdots \\
                  z_{2x}=z_{(2x-1)1}-z_{(2x-2)2}+\cdots+z_{1(2x-1)} etc.
                \end{array} \right.  \label{zrelations}
\end{E}where $\alpha, \beta$ are compositions.
Note that the last equation in \eqref{zrelations} is equivalent to
\begin{E} z_{xx}=(-1)^{x+1}(z_{2x}-z_{(2x-1)1}+z_{(2x-2)2}-\cdots\widehat{z_{xx}}\cdots-z_{1(2x-1)})  \label{zxx}.\end{E}
Let $\gamma$ be a composition with length $k$. Using the first equation in \eqref{zrelations}, we have 
\begin{E} z_{\alpha\cdot\beta}+z_{\alpha\odot\beta}=z_{\beta\cdot\alpha}+z_{\beta\odot\alpha} \label{switch}.\end{E}
By \cite[Proposition 2.2]{HDL} we can sort $z_\gamma$, that is $z_\gamma=z_{\lambda(\gamma)}+$ a sum of $z_\delta$'s with $\delta$ having $k-1$ or fewer parts. For $\alpha=\alpha_1\cdots\alpha_m$, define $prod(\alpha)$ to be the product of the parts of the composition $\alpha$, that is $prod(\alpha)=\alpha_1\times\alpha_2\times\cdots\times\alpha_m$. The partition $\alpha$ is called a {\it semi-strict} partition if it can be written in the form $\alpha=\alpha_1\alpha_2\cdots\alpha_k1\cdots1$ where $\alpha_1\alpha_2\cdots\alpha_k$ is a strict partition. 

Suppose that $\lambda(\gamma)=g_1g_2\ldots g_k$. If there is no $i$, $1\leq i\leq k-1$, such that $g_i=g_{i+1}=t>1$ then $\lambda(\gamma)$ is a semi-strict partition and we have \eqref{zgamma}, otherwise
\begin{E}\begin{array}{lll}
 z_\gamma &=&z_{\lambda(\gamma)}+\hbox{\scriptsize a sum of $z_\delta$'s such that $\ell(\delta)<k$} \\
 & = & z_{g_ig_{i+1}\ldots g_kg_1\ldots g_{i-1}}+\hbox{\scriptsize a sum of $z_\delta$'s such that $\ell(\delta)<k$} \\
	& =& z_{g_ig_{i+1}}z_{g_{i+2}\ldots g_kg_1\ldots g_{i-1}}+	\hbox{\scriptsize a sum of $z_\delta$'s such that $\ell(\delta)<k$} \\
	& =& (-1)^{t+1}[(z_{2t}-z_{(2t-1)1}+z_{(2t-2)2}-\cdots\widehat{z_{tt}}\cdots-z_{1(2t-1)})z_{g_{i+2}\ldots g_kg_1\ldots g_{i-1}}]\\
	&&+	\hbox{\scriptsize a sum of $z_\delta$'s such that $\ell(\delta)<k$}\\
	&= & (-1)^{t+1}[-z_{(2t-1)1g_{i+2}\ldots g_kg_1\ldots g_{i-1}}+z_{(2t-2)2g_{i+2}\ldots g_kg_1\ldots g_{i-1}}\\
	&&-\cdots \widehat{z_{ttg_{i+2}\ldots g_kg_1\cdots g_{i-1}}}
	\cdots- z_{1(2t-1)g_{i+2}\ldots g_kg_1\ldots g_{i-1}}]+ 	\hbox{\scriptsize a sum of $z_\delta$'s such that $\ell(\delta)<k$}\\
	&=& (-1)^{t+1}[-z_{\lambda((2t-1)1g_{i+2}\ldots g_kg_1\ldots g_{i-1})}+z_{\lambda((2t-2)2g_{i+2}\ldots g_kg_1\ldots g_{i-1})}\\
	&&-\cdots\widehat{z_{\lambda(ttg_{i+2}\ldots g_kg_1\ldots g_{i-1})}}
	\cdots- z_{\lambda(1(2t-1)g_{i+2}\ldots g_kg_1\ldots g_{i-1})}]+	\hbox{\scriptsize a sum of $z_\delta$'s such that $\ell(\delta)<k$}\\
\end{array}\label{process}\end{E}where we used  \eqref{switch} for the second, the first equation of \eqref{zrelations} for the third, \eqref{zxx} for the fourth, the first equation of \eqref{zrelations} for the fifth, and sorting for the sixth equality. Although 
$\lambda((2t-1)1g_{i+2}\ldots g_kg_1\ldots g_{i-1})$, $\lambda((2t-2)2g_{i+2}\ldots g_kg_1\ldots g_{i-1})$, $\ldots$, $\lambda(1(2t-1)g_{i+2}\ldots g_kg_1\ldots g_{i-1})$ have $k$ parts, the product of their parts is smaller than $prod(\lambda(\gamma))$ since $(2t-1),2\times(2t-2),\ldots,2t-1<t^2$. 
We repeat the process in \eqref{process} for each of the terms with $k$ parts in the last line of \eqref{process}. Since $prod(\alpha)$ is  a positive integer,   the process terminates, which yields 
\begin{E} z_\gamma= (\hbox{\scriptsize a sum of $z_\sigma$'s such that $\sigma$ is a semi-strict partition with $\ell(\sigma)=k$})\hspace{5pt}+\hspace{5pt}(\hbox{\scriptsize a sum of $z_\delta$'s such that $\ell(\delta)<k$})  \label{zgamma}.\end{E}
Now if $\sigma$ is a semi-strict partition with at least two 1's, that is $\sigma=\sigma'11$ where $\sigma'$ is a semi-strict partition and $\ell(\sigma')=k-2$, then we have
\begin{E}\begin{array}{lll}
z_\sigma & =&z_{11\sigma'}+\hbox{\scriptsize a sum of $z_\delta$'s such that $\ell(\delta)<k$} \\
 &=& z_{11}z_{\sigma'}+\hbox{\scriptsize a sum of $z_\delta$'s such that $\ell(\delta)<k$}\\
 &=&z_2z_{\sigma'}+\hbox{\scriptsize a sum of $z_\delta$'s such that $\ell(\delta)<k$} \\
 &=&z_{2\sigma'}+z_{2\odot\sigma'}+\hbox{\scriptsize a sum of $z_\delta$'s such that $\ell(\delta)<k$}
 \end{array} \label{zsigma} \end{E}where we used \eqref{switch} for the first, the first equation of \eqref{zrelations} for the second, the second equation of \eqref{zrelations} for the third, and the first equation of \eqref{zrelations} for the fourth equality.  
Note that $\ell(2\sigma')=k-1$ and $\ell(2\odot\sigma')=k-2$. If $\sigma$ does not have two 1's then it is a strict partition. Now applying \eqref{zsigma} to each  $z_\sigma$ with $\sigma $ having at least two 1's in \eqref{zgamma}, we have 
$$ z_\gamma= (\hbox{\scriptsize a sum of $z_\sigma$'s such that $\sigma$ is a strict partition with $\ell(\sigma)=k$})\hspace{5pt}+\hspace{5pt} (\hbox{\scriptsize a sum of $z_\delta$'s such that $\ell(\delta)<k$}).$$
A trivial induction on the length of $\gamma$ now shows that any $z_\gamma$ in the quotient can be written as a linear combination of $z_\lambda$, $\lambda \vdash n$ and $\lambda$ is a strict partition. \end{proof}

However, this is not the only possible set of relations and we now develop another set. This alternative set will help simplify some of our subsequent proofs in addition to being of independent interest.

\begin{theorem} 
\label{ribbonrelations2} 
Let {$z_{\alpha}, \alpha \vDash n, n\ge 1$} be commuting indeterminates.  Then 
as algebras, ${\Omega}$ is isomorphic to the quotient 
$${\Q[z_{\alpha}]/ 
 \langle z_{\alpha}\ z_{\beta}-z_{\alpha\cdot \beta} - z_{\alpha \odot \beta}, \xi _2, \xi_4, \ldots\rangle}$$where
 $\xi_{2m}$ is the even transpose form $\xi _{2m} = z_{2m} - z_{\underbrace{1\ldots1}_{2m}}$. Thus, all relations amongst ribbon Schur $Q$-functions are generated by $\Qr_{\alpha}\ \Qr_{\beta}= \Qr_{\alpha\cdot \beta} + \Qr_{\alpha \odot \beta}$ and  $\Qr_{2m} = \Qr_{\underbrace{1\ldots1}_{2m}}$, $m\geq 1$.

\end{theorem} 

We devote the next subsection to the proof of this theorem.

\subsection{Equivalence of relations}

We say that the set of relationships $A$ {\it implies} the set of relationships $B$, if we can deduce $B$ from $A$. Two sets of relationships are {\it equivalent}, if each one implies the other one. 
%Let ${\mathfrak C}_n$ denote the set of all coarsenings of $\underbrace{1\ldots1}_{n}$, that is, the set of all compositions of $n$. For a composition $\beta$, let $\beta'$ be the composition that satisfies $\beta=\beta_1\cdot\beta'$; in other words, if $\beta=\beta_1\beta_2\ldots\beta_k$ then $\beta'=\beta_2\ldots\beta_k$.
\begin{itemize}
\item For all compositions $\alpha$ and $\beta$, refer to 
$$z_\alpha z_\beta=z_{\alpha\cdot\beta}+z_{\alpha\odot\beta}$$
as multiplication. 
\item For all positive integers $x$, refer to the set of 
$$z_{2x}=z_{(2x-1)1}-z_{(2x-2)2}+\cdots-z_{2(2x-2)}+z_{1(2x-1)}$$
as $EE$.
\item For all positive integers $x$, refer to the set of 
$$2z_{2x}=z_{2x-1}z_1-z_{2x-2}z_2+\cdots-z_2z_{2x-2}+z_1z_{2x-1}$$
as $EI$.
\item  For all positive integers $x$, refer to the set of 
$$z_{x}=z_{\underbrace{1\ldots1}_{x}}$$
as $T$. 
\item  For all positive integers $x$, refer to the set of 
$$z_{2x}=z_{\underbrace{1\ldots1}_{2x}}$$
as $ET$. 
\end{itemize}

\begin{lemma}
Multiplication and $EE$ is equivalent to multiplication and $EI$.
\label{eeei} \end{lemma}

\begin{proof} 
$$\begin{array}{ll}
                &z_{2x}=z_{(2x-1)1}-z_{(2x-2)2}+\cdots-z_{2(2x-2)}+z_{1(2x-1)}\\
\Leftrightarrow &z_{2x}=(z_{2x-1}z_1-z_{2x})-(z_{2x-2}z_2-z_{2x})+\cdots-(z_2z_{2x-2}-z_{2x})+(z_1z_{2x-1}-z_{2x})\\
\Leftrightarrow & 2z_{2x}=z_{2x-1}z_1-z_{2x-2}z_2+\cdots-z_2z_{2x-2}+z_1z_{2x-1} 
\end{array}$$
where we used multiplication for the first  equivalence.\end{proof}

\begin{lemma}
Multiplication and $T$ is equivalent to multiplication and $EI$.
\label{tei}\end{lemma}

\begin{proof} First we show that the set of $T$ and multiplication implies $EI$.
$$
\begin{array}{ll}
  & z_{2x-1}z_1-z_{2x-2}z_2+z_{2x-3}z_3-\cdots-z_2z_{2x-2}+z_1z_{2x-1}\\
= & z_{2x-1}z_1-z_{2x-2}z_{11}+z_{2x-3}z_{111}-\cdots-z_2z_{\underbrace{1\ldots1}_{2x-2}}+z_1z_{\underbrace{1\ldots1}_{2x-1}}\\
= & (z_{2x}+z_{(2x-1)1})-(z_{(2x-1)1}+z_{(2x-2)11})+(z_{(2x-2)11}+z_{(2x-3)111})-\cdots-\\
 &(z_{3\underbrace{1\ldots1}_{2x-3}}+z_{2\underbrace{1\ldots1}_{2x-2}})+(z_{2\underbrace{1\ldots1}_{2x-2}}+z_{\underbrace{1\ldots1}_{2x}})\\
= & z_{2x}+ z_{\underbrace{1\ldots1}_{2x}}\\
= & 2z_{2x}
\end{array}
$$
where we used $T$ for the first, multiplication for the second, and $T$ for the fourth equality. 

Now we proceed by induction to show that the set of $EI$ and multiplication implies $T$. The base case is $z_1=z_1$.
% First we check the base case
%$$2z_2=z_1z_1=z_2+z_{11}\Longrightarrow z_2=z_{11}.$$
Assume the assertion is true for all $n$ smaller than $2x$, so the set of $EI$ and multiplication implies $z_n=z_{\underbrace{1\ldots1}_{n}}$ for all $n<2x$. We show that it is true for $2x$ and $2x+1$ as well. 
$$ \begin{array}{lll}
2z_{2x} & = & z_{2x-1}z_1-z_{2x-2}z_2+z_{2x-3}z_3-\cdots-z_2z_{2x-2}+z_1z_{2x-1}\\
        & = & z_{2x-1}z_1-z_{2x-2}z_{11}+z_{2x-3}z_{111}-\cdots-z_2z_{\underbrace{1\ldots1}_{2x-2}}+z_1z_{\underbrace{1\ldots1}_{2x-1}}\\
        & = & (z_{2x}+z_{(2x-1)1})-(z_{(2x-1)1}+z_{(2x-2)11})+(z_{(2x-2)11}+z_{(2x-3)111})-\cdots-\\
&      &(z_{3\underbrace{1\ldots1}_{2x-3}}+z_{2\underbrace{1\ldots1}_{2x-2}})+(z_{2\underbrace{1\ldots1}_{2x-2}}+z_{\underbrace{1\ldots1}_{2x}})\\
				& = & z_{2x}+ z_{\underbrace{1\ldots1}_{2x}}\\
\end{array}$$
where we used $EI$ for the first, the induction hypothesis for the second, and multiplication for the third equality. Thus $z_{2x}=z_{\underbrace{1\ldots1}_{2x}}$. Now we show that $z_{2x+1}=z_{\underbrace{1\ldots1}_{2x+1}}$.
$$\begin{array}{lll}
0 & =& z_{2x}z_1-z_{2x-1}z_2+z_{2x-2}z_3-\cdots+z_2z_{2x-1}-z_1z_{2x}\\
  & =& z_{2x}z_1-z_{2x-1}z_{11}+z_{2x-2}z_{111}-\cdots+z_2z_{\underbrace{1\ldots1}_{2x-1}}-z_1z_{\underbrace{1\ldots1}_{2x}}\\
  & =& (z_{2x+1}+z_{(2x)1})-(z_{(2x)1}+z_{(2x-1)11})+(z_{(2x-1)11}+z_{(2x-2)111})-\cdots+\\
  &  & (z_{3\underbrace{1\ldots1}_{2x-2}}+z_{2\underbrace{1\ldots1}_{2x-1}})-(z_{2\underbrace{1\ldots1}_{2x-1}}+z_{\underbrace{1\ldots1}_{2x+1}})\\
  & =& z_{2x+1}-z_{\underbrace{1\ldots1}_{2x+1}}
\end{array}$$ 
where we used the induction hypothesis and $z_{2x}=z_{\underbrace{1\ldots1}_{2x}}$ for the second, and multiplication for the third equality. Thus $z_{2x+1}=z_{\underbrace{1\ldots1}_{2x+1}}$, which completes the induction. \end{proof}

\begin{lemma}
Multiplication and $T$ is equivalent to multiplication and $ET$.
\label{tet}\end{lemma}

\begin{proof} The set of relationships $ET$ is a subset of $T$, thus $T$ implies $ET$. To prove the converse, we need to show $z_{2x+1}=z_{\underbrace{1\ldots1}_{2x+1}}$ given $ET$ and multiplication. We proceed by induction. The base case is   $z_1=z_1$. Assume the result is true for all odd positive integers smaller than $2x+1$, then

$$\begin{array}{lll}
0 & =& z_{2x}z_1-z_{2x-1}z_2+z_{2x-2}z_3-\cdots+z_2z_{2x-1}-z_1z_{2x}\\
  & =& z_{2x}z_1-z_{2x-1}z_{11}+z_{2x-2}z_{111}-\cdots+z_2z_{\underbrace{1\ldots1}_{2x-1}}-z_1z_{\underbrace{1\ldots1}_{2x}}\\
  & =& (z_{2x+1}+z_{(2x)1})-(z_{(2x)1}+z_{(2x-1)11})+(z_{(2x-1)11}+z_{(2x-2)111})-\cdots+\\
  &  & (z_{3\underbrace{1\ldots1}_{2x-2}}+z_{2\underbrace{1\ldots1}_{2x-1}})-(z_{2\underbrace{1\ldots1}_{2x-1}}+z_{\underbrace{1\ldots1}_{2x+1}})\\
  & =& z_{2x+1}-z_{\underbrace{1\ldots1}_{2x+1}}
\end{array}$$
where we used $ET$ and the induction hypothesis for the second, and multiplication for the third equality. Thus $z_{2x+1}=z_{\underbrace{1\ldots1}_{2x+1}}$, which completes the induction. \end{proof}

Combining Lemma \ref{eeei}, Lemma \ref{tei} and Lemma \ref{tet} we get

\begin{proposition}
Multiplication and $EE$ is equivalent to multiplication and $ET$.
\label{coreeet}\end{proposition}

Theorem~\ref{ribbonrelations2} now follows from Theorem~\ref{ribbonrelations} and Proposition~\ref{coreeet}.

%%%%%% Equality Section %%%%%%%%%%
\section{\texorpdfstring{Equality of ordinary skew Schur $Q$-functions}{Equality of ordinary skew Schur Q-functions}}\label{sec:eqskewschurq}

We now turn our attention to determining when two ordinary skew Schur $Q$-functions are equal. Illustrative examples of the results in this section can be found in the next section, when we restrict our attention to ribbon Schur $Q$-functions. In order to prove our main result on equality, Theorem~\ref{the:bigone}, which is analogous to \cite[Theorem 7.6]{HDL2}, we need to prove an analogy of \cite[Proposition 2.1]{HDL}. First we need to prove a  Jacobi-Trudi style determinant formula.

Let $D_1,D_2,\ldots,D_k$ denote skew diagrams, and recall from Section~\ref{sec:diagrams} that $$D_1\bigstar_1D_2\bigstar_2D_3\bigstar_3\cdots\bigstar_{k-1}D_k$$in which $\bigstar_i$ is either $\cdot$ or $\odot$ is a well-defined skew diagram. Set
$$\bar{\bigstar}_i=\left \{
\begin{array}{ll}
\odot & \hbox{if }\bigstar_i=\cdot\\
\cdot & \hbox{if }\bigstar_i=\odot .
\end{array} \right. $$
With this in mind we have

\begin{proposition} Let $s_D$ denote the skew Schur function indexed by the ordinary skew diagram $D$. Then
$$s_{D_1\bigstar_1D_2\bigstar_2D_3\bigstar_3\cdots\bigstar_{k-1}D_k}=
\det \left [ \begin{array}{ccccc}
                                      s_{D_1} & s_{D_1\bar{\bigstar}_1D_2}		& s_{D_1\bar{\bigstar}_1D_2\bar{\bigstar}_2D_3}	& \cdots	& s_{D_1\bar{\bigstar}_1D_2\bar{\bigstar}_2\cdots\bar{\bigstar}_{k-1}D_k}\\ 
                                      1    &	s_{D_2}	&	s_{D_2\bar{\bigstar}_2D_3}	&	\cdots	&	s_{D_2\bar{\bigstar}_2\bar{\bigstar}_3\cdots\bar{\bigstar}_{k-1}D_k} \\
                                      			&	1			&	s_{D_3}	&	\cdots	& s_{D_3\bar{\bigstar}_3\cdots\bar{\bigstar}_{k-1}D_k} \\
                                      			&				& \ddots 	& 			&	\vdots \\
                                       0    &				&					&1			&s_{D_k} \end{array} \right ]  		.$$
\label{propsd}\end{proposition}

\begin{proof} We proceed by induction on $k$. Assuming the assertion is true for $k-1$, we show that it is true for $k$ as well. Note that the base case, $k=2$, is, say \cite[Proposition 4.1]{HDL2}, that 
\begin{E} s_{D_1}s_{D_2}=s_{D_1\cdot D_2}+s_{D_1\odot D_2} \label{skewschur0}\end{E}for skew diagrams $D_1, D_2$.

By the induction hypothesis, we have 
\begin{E}
\det \left [ \begin{array}{cccc}
                                      s_{D} & s_{D\bar{\bigstar}_2D_3}		& \cdots	& s_{D\bar{\bigstar}_2D_3\bar{\bigstar}_3\cdots\bar{\bigstar}_{k-1}D_k}\\ 
                                      1    &	s_{D_3}	                    &	\cdots	&	
s_{D_3\bar{\bigstar}_3\cdots\bar{\bigstar}_{k-1}D_k} \\
                                      			&	\ddots			&		&	\vdots \\
                                       0    &									&1			&s_{D_k} \end{array} \right ]  	=s_{D\bigstar_2D_3\bigstar_3\cdots\bigstar_{k-1}D_k}
\label{skewschur1}\end{E}
where $D$ can be any  skew diagram. Now expanding over the first column yields 
\begin{E}\begin{array}{ll}
\det \left [ \begin{array}{ccccc}
                                      s_{D_1} & s_{D_1\bar{\bigstar}_1D_2}		& s_{D_1\bar{\bigstar}_1D_2\bar{\bigstar}_2D_3}	& \cdots	& s_{D_1\bar{\bigstar}_1D_2\bar{\bigstar}_2\cdots\bar{\bigstar}_{k-1}D_k}\\ 
                                      1    &	s_{D_2}	&	s_{D_2\bar{\bigstar}_2D_3}	&	\cdots	&	s_{D_2\bar{\bigstar}_2\bar{\bigstar}_3\cdots\bar{\bigstar}_{k-1}D_k} \\
                                      			&	1			&	s_{D_3}	&	\cdots	& s_{D_3\bar{\bigstar}_3\cdots\bar{\bigstar}_{k-1}D_k} \\
                                      			&				& \ddots 	& 			&	\vdots \\
                                       0    &				&					&1			&s_{D_k} \end{array} \right ] &=\\ %%
                                      s_{D_1}\times \det \left [ \begin{array}{cccc}
                                      s_{D_2} & s_{D_2\bar{\bigstar}_2D_3}		& \cdots	& s_{D_2\bar{\bigstar}_2D_3\bar{\bigstar}_3\cdots\bar{\bigstar}_{k-1}D_k}\\ 
                                      1    &	s_{D_3}	                    &	\cdots	&	
s_{D_3\bar{\bigstar}_3\cdots\bar{\bigstar}_{k-1}D_k} \\
                                      			&	\ddots			&		&	\vdots \\
                                       0    &									&1			&s_{D_k} \end{array} \right ]  & -\\ %% 
                                      \det \left [ \begin{array}{cccc}
                                      s_{D_1\bar{\bigstar}_1D_2} & s_{D_1\bar{\bigstar}_1D_2\bar{\bigstar}_2D_3}		& \cdots	& s_{D_1\bar{\bigstar}_1D_2\bar{\bigstar}_2D_3\bar{\bigstar}_3\cdots\bar{\bigstar}_{k-1}D_k}\\ 
                                      1    &	s_{D_3}	                    &	\cdots	&	
s_{D_3\bar{\bigstar}_3\cdots\bar{\bigstar}_{k-1}D_k} \\
                                      			&	\ddots			&		&	\vdots \\
                                       0    &									&1			&s_{D_k} \end{array} \right ]  \ .
\label{skewschur2}\end{array}\end{E}\\

Note that the first and second determinant on the right side of \eqref{skewschur2} are equal to the determinant in \eqref{skewschur1} for, respectively, $D=D_2$ and $D=D_1\bar{\bigstar}_1D_2$. Thus, the equality in \eqref{skewschur1} implies that \eqref{skewschur2} is equal to 
$$s_{D_1}\times s_{D_2\bigstar_2D_3\bigstar_3\cdots\bigstar_{k-1}D_k}-s_{D_1\bar{\bigstar}_1D_2\bigstar_2D_3\bigstar_3\cdots\bigstar_{k-1}D_k}$$and
because of \eqref{skewschur0}, the last expression is equal to 
$$s_{D_1\bigstar_1D_2\bigstar_2D_3\bigstar_3\cdots\bigstar_{k-1}D_k}.$$
This completes the induction.\end{proof}

Let $\alpha$ be a ribbon such that $$\alpha =\square\bigstar _1 \square \bigstar _2 \cdots \bigstar _{k-1} \square$$and $|\alpha|=k$. In Proposition~\ref{propsd} set $D_i=D$ for $i$ odd and $D_i=D^t$ for $i$ even for some skew diagram $D$ so that $D\bigstar_1 D^t \bigstar_2 D\bigstar_3\cdots = \alpha\bullet D$. Note that 
$$\alpha^t\bullet D=D\bar{\bigstar}_{k-1}D^t\bar{\bigstar}_{k-2}D\bar{\bigstar}_{k-3}\cdots$$
therefore, 
$$(\alpha^t)^\circ\bullet D=D\bar{\bigstar}_1D^t\bar{\bigstar}_2D\bar{\bigstar}_3\cdots.$$
Using Proposition~\ref{propsd} with the above setting, we have the following corollary.

\begin{corollary}
$$s_{\alpha \bullet D}=\det \left [ \begin{array}{ccccc}
                                      \ast & \ast		& \ast	& \cdots	& s_{(\alpha^t)^\circ\bullet D}\\ 
                                      1    &	\ast	&	\ast	&	\cdots	&	\ast \\
                                      			&	1			&	\ast	&	\cdots	& \ast \\
                                      			&				& \ddots 	& 			&	\vdots \\
                                       0    &				&					&1			&\ast \end{array} \right ]  			  $$where the skew Schur functions indexed by skew diagrams with fewer than $|\alpha|$ blocks of $D$ or $D^t$ are denoted by $\ast$. 
\label{dihamel} \end{corollary}

We are now ready to derive our first ordinary skew Schur $Q$-function equalities.

\begin{proposition} If $\alpha$ is a ribbon and $D$ is a skew diagram then ${\mathfrak s}_{\alpha \bullet  D }={\mathfrak s}_{\alpha^\circ \bullet  D }$ and ${\mathfrak s}_{\alpha \bullet  D }={\mathfrak s}_{\alpha \bullet  D ^t}.$ \label{diproprottrans}\end{proposition}

\begin{proof} We induct on $|\alpha|$. The base case is easy as $1=1^\circ$ and ${\mathfrak s}_ D ={\mathfrak s}_{ D ^t}$ by \eqref{eq:Qtr}. Assume the proposition is true for $|\alpha|<n$. We first show that ${\mathfrak s}_{\alpha \bullet  D }={\mathfrak s}_{\alpha^\circ \bullet  D }$ for all $\alpha$'s with $|\alpha|=n$, by inducting on the number of parts in $\alpha$, that is $\ell(\alpha)$. The base case, $\ell(\alpha)=1$, is straightforward as $n=n^\circ$. Assume ${\mathfrak s}_{\alpha \bullet  D }={\mathfrak s}_{\alpha^\circ \bullet  D }$ is true for all compositions $\alpha$ with fewer than $k$ parts (the hypothesis for the second induction). Let  $\alpha=\alpha_1\cdots\alpha_k$.  
Using \eqref{eq:Qmult}, we know that for all skew diagrams $V$ and $L$ we have 
$${\mathfrak s}_{ V\cdot L }={\mathfrak s}_ V{\mathfrak s}_L -{\mathfrak s}_{ V\odot L } . $$
We consider the following four cases. Note that in each case we set $V$ and $L$ such that  $ V\cdot L =\alpha_1\cdots \alpha_{k-1}\alpha_k \bullet D =\alpha\bullet D$ and $V\odot L =\alpha_1\cdots(\alpha_{k-1}+\alpha_k) \bullet D$. Also, note that since $|\alpha_1\cdots\alpha_{k-1}|<n$ and $|\alpha_k|<n$, we can use the induction hypothesis of the first induction (i.e. we can rotate the first and transpose the second component). Furthermore, even though $|\alpha_1\cdots(\alpha_{k-1}+\alpha_k)|=n$, the number of parts in $\alpha_1\cdots(\alpha_{k-1}+\alpha_k)$ is $k-1$ and therefore we can use the induction hypothesis of the second induction:
\\
{\it Case 1: $|\alpha_1\cdots \alpha_{k-1}|$ is even and $|\alpha_k|$ is even.} Set $V=\alpha_1 \cdots \alpha_{k-1} \bullet D$ and $L=\alpha_k \bullet D$. Then 
\begin{eqnarray*}{\mathfrak s}_{\alpha \bullet D}&=&{\mathfrak s}_{\alpha_1\cdots\alpha_{k-1}\bullet  D } {\mathfrak s}_{\alpha_k \bullet  D }-{\mathfrak s}_{\alpha_1\cdots(\alpha_{k-1}+\alpha_k)\bullet  D }\\&=&{\mathfrak s}_{\alpha_k \bullet  D }{\mathfrak s}_{\alpha_{k-1}\cdots\alpha_{1}\bullet  D }-{\mathfrak s}_{(\alpha_k+\alpha_{k-1})\cdots \alpha_1\bullet  D }={\mathfrak s}_{\alpha_k\alpha_{k-1}\cdots\alpha_{1}\bullet  D }={\mathfrak s}_{\alpha^\circ \bullet  D }. \end{eqnarray*}
\\
{\it Case 2: $|\alpha_1\cdots \alpha_{k-1}|$ is even and $|\alpha_k|$ is odd.} Set $ V =\alpha_1 \cdots \alpha_{k-1} \bullet  D $ and $ L =\alpha_k \bullet  D $. Then 
\begin{eqnarray*}{\mathfrak s}_{\alpha \bullet  D }&=&{\mathfrak s}_{\alpha_1\cdots\alpha_{k-1}\bullet  D } {\mathfrak s}_{\alpha_k \bullet  D }-{\mathfrak s}_{\alpha_1\cdots(\alpha_{k-1}+\alpha_k)\bullet  D }\\&=&{\mathfrak s}_{\alpha_k \bullet  D }{\mathfrak s}_{\alpha_{k-1}\cdots\alpha_{1}\bullet  D ^t} -{\mathfrak s}_{(\alpha_k+\alpha_{k-1})\cdots \alpha_1\bullet  D }={\mathfrak s}_{\alpha_k\alpha_{k-1}\cdots\alpha_{1}\bullet  D }={\mathfrak s}_{\alpha^\circ \bullet  D }. \end{eqnarray*}
\\
{\it Case 3: $|\alpha_1\cdots \alpha_{k-1}|$ is odd and $|\alpha_k|$ is even.} Set $ V =\alpha_1 \cdots \alpha_{k-1} \bullet  D $ and $ L =\alpha_k \bullet  D ^t$. Then  
\begin{eqnarray*}{\mathfrak s}_{\alpha \bullet  D }&=&{\mathfrak s}_{\alpha_1\cdots\alpha_{k-1}\bullet  D } {\mathfrak s}_{\alpha_k \bullet  D ^t}-{\mathfrak s}_{\alpha_1\cdots(\alpha_{k-1}+\alpha_k)\bullet  D }\\&=&{\mathfrak s}_{\alpha_k \bullet  D }{\mathfrak s}_{\alpha_{k-1}\cdots\alpha_{1}\bullet  D }-{\mathfrak s}_{(\alpha_k+\alpha_{k-1})\cdots \alpha_1\bullet  D }={\mathfrak s}_{\alpha_k\alpha_{k-1}\cdots\alpha_{1}\bullet  D }={\mathfrak s}_{\alpha^\circ \bullet  D }. \end{eqnarray*}
\\
{\it Case 4: $|\alpha_1\cdots \alpha_{k-1}|$ is odd and $|\alpha_k|$ is odd.} Set $ V =\alpha_1 \cdots \alpha_{k-1} \bullet  D $ and $ L =\alpha_k \bullet  D ^t$. Then  
\begin{eqnarray*}{\mathfrak s}_{\alpha \bullet  D }&=&{\mathfrak s}_{\alpha_1\cdots\alpha_{k-1}\bullet  D } {\mathfrak s}_{\alpha_k \bullet  D ^t}-{\mathfrak s}_{\alpha_1\cdots(\alpha_{k-1}+\alpha_k)\bullet  D }\\&=&{\mathfrak s}_{\alpha_k \bullet  D }{\mathfrak s}_{\alpha_{k-1}\cdots\alpha_{1}\bullet  D ^t}-{\mathfrak s}_{(\alpha_k+\alpha_{k-1})\cdots \alpha_1\bullet  D }={\mathfrak s}_{\alpha_k\alpha_{k-1}\cdots\alpha_{1}\bullet  D }={\mathfrak s}_{\alpha^\circ \bullet  D }. \end{eqnarray*}

This completes the second induction. Now to complete the first induction, we show that ${\mathfrak s}_{\alpha \bullet  D }={\mathfrak s}_{\alpha \bullet  D ^t}$ where $|\alpha|=n$. 
\\
Suppose $n$ is odd. By Corollary \ref{dihamel}, we have 
$$s_{\alpha \bullet  D }=\det \left [ \begin{array}{ccccc}
                                      \ast & \ast		& \ast	& \cdots	& s_{(\alpha^t)^\circ\bullet  D }\\ 
                                      1    &	\ast	&	\ast	&	\cdots	&	\ast \\
                                      			&	1			&	\ast	&	\cdots	& \ast \\
                                      			&				& \ddots 	& 			&	\vdots \\
                                      0     &				&					&1			&\ast \end{array} \right ]  	.		  $$                
Expanding the above determinant we have 
$$ s_{\alpha \bullet  D }=X+s_{(\alpha^t)^\circ\bullet  D } $$
where $X$ is comprised of skew Schur functions indexed by skew diagrams with fewer than $|\alpha|$ blocks of $ D $ or $ D ^t$. Applying $\theta$ to both sides of the above equation yields 
\begin{E} {\mathfrak s}_{\alpha \bullet  D }={\mathfrak X}+{\mathfrak s}_{(\alpha^t)^\circ\bullet  D }={\mathfrak X}+{\mathfrak s}_{\alpha^t\bullet  D }= {\mathfrak X}+{\mathfrak s}_{(\alpha^t\bullet  D )^t}={\mathfrak X}+{\mathfrak s}_{\alpha\bullet  D ^t} \label{ditransodd1} 
\end{E}where we used the result of the second induction for the second, \eqref{eq:Qtr} for the third and \eqref{ditransposition} for the fourth equality.
 
Similarly, 
$$s_{\alpha \bullet  D ^t}=\det \left [ \begin{array}{ccccc}
                                      \ast & \ast		& \ast	& \cdots	& s_{(\alpha^t)^\circ\bullet  D ^t}\\ 
                                      1    &	\ast	&	\ast	&	\cdots	&	\ast \\
                                      			&	1			&	\ast	&	\cdots	& \ast \\
                                      			&				& \ddots 	& 			&	\vdots \\
                                      0     &				&					&1			&\ast \end{array} \right ]  			  $$and expanding the determinant we have 
$$ s_{\alpha \bullet  D ^t}=X'+s_{(\alpha^t)^\circ\bullet  D ^t} 
$$where $X'$ is again comprised of skew Schur functions indexed by skew diagrams with fewer than $|\alpha|$ blocks of $ D $ or $ D ^t$. By the induction hypothesis of the first induction (i.e. the induction on $|\alpha|$), we can assume $\theta(X')=\theta(X)={\mathfrak X}$. Now we apply $\theta$ to both sides of the above equation, thus
\begin{E} {\mathfrak s}_{\alpha \bullet  D ^t}={\mathfrak X}+{\mathfrak s}_{(\alpha^t)^\circ\bullet  D ^t}={\mathfrak X}+{\mathfrak s}_{\alpha^t\bullet  D ^t}= {\mathfrak X}+{\mathfrak s}_{(\alpha^t\bullet  D ^t)^t}={\mathfrak X}+{\mathfrak s}_{\alpha\bullet  D } \label{ditransodd2} \end{E}where, again, we used the result of the second induction for the second, \eqref{eq:Qtr} for the third and \eqref{ditransposition} for the fourth equality. Now \eqref{ditransodd1} and \eqref{ditransodd2} imply ${\mathfrak s}_{\alpha\bullet  D }={\mathfrak s}_{\alpha \bullet  D ^t}$ for the case $|\alpha|=n$ odd.
 
\comment{Now suppose $n$ is even. 
Again by Corollary \ref{dihamel}, we have 
$$s_{\alpha \bullet  D }=\det \left [ \begin{array}{ccccc}
                                      \ast & \ast		& \ast	& \cdots	& s_{(\alpha^t)^\circ\bullet  D }\\ 
                                      1    &	\ast	&	\ast	&	\cdots	&	\ast \\
                                      			&	1			&	\ast	&	\cdots	& \ast \\
                                      			&				& \ddots 	& 			&	\vdots \\
                                     0      &				&					&1			&\ast \end{array} \right ]  		.	  $$                
Expanding the above determinant we have 
$$ s_{\alpha \bullet  D }=Y-s_{(\alpha^t)^\circ\bullet  D } $$
where $Y$ is comprised of skew Schur functions indexed by skew diagrams with fewer than $|\alpha|$ blocks of $ D $ or $ D ^t$. Applying $\theta$ to both sides of the above equation yields 
\begin{E} {\mathfrak s}_{\alpha \bullet  D }={\mathfrak Y}-{\mathfrak s}_{(\alpha^t)^\circ\bullet  D }={\mathfrak Y}-{\mathfrak s}_{\alpha^t\bullet  D }= {\mathfrak Y}-{\mathfrak s}_{(\alpha^t\bullet  D )^t}={\mathfrak Y}-{\mathfrak s}_{\alpha\bullet  D } \label{ditranseven1} \end{E}where we used the result of the second induction for the second, \eqref{eq:Qtr} for the third and \eqref{ditransposition} for the fourth equality.
 
Similarly, 
$$s_{\alpha \bullet  D ^t}=\det \left [ \begin{array}{ccccc}
                                      \ast & \ast		& \ast	& \cdots	& s_{(\alpha^t)^\circ\bullet  D ^t}\\ 
                                      1    &	\ast	&	\ast	&	\cdots	&	\ast \\
                                      			&	1			&	\ast	&	\cdots	& \ast \\
                                      			&				& \ddots 	& 			&	\vdots \\
                                      0     &				&					&1			&\ast \end{array} \right ]  			  $$and expanding the determinant we have 
$$ s_{\alpha \bullet  D ^t}=Y'-s_{(\alpha^t)^\circ\bullet  D ^t} 
$$where $Y'$ is again comprised of skew Schur functions indexed by skew diagrams with fewer than $|\alpha|$ blocks of $ D $ or $ D ^t$. By the induction hypothesis on $|\alpha|$ we can assume $\theta(Y')=\theta(Y)={\mathfrak Y}$. Now we apply $\theta$ to both sides of the above equation, thus
\begin{E} {\mathfrak s}_{\alpha \bullet  D ^t}={\mathfrak Y}-{\mathfrak s}_{(\alpha^t)^\circ\bullet  D ^t}={\mathfrak Y}-{\mathfrak s}_{\alpha^t\bullet  D ^t}= {\mathfrak Y}-{\mathfrak s}_{(\alpha^t\bullet  D ^t)^t}={\mathfrak Y}-{\mathfrak s}_{\alpha\bullet  D ^t} \label{ditranseven2} 
\end{E}where, again, we used the result of the second induction for the second, \eqref{eq:Qtr} for the third and \eqref{ditransposition} for the fourth equality. Now \eqref{ditranseven1} and \eqref{ditranseven2} imply ${\mathfrak s}_{\alpha\bullet  D }={\mathfrak s}_{\alpha \bullet  D ^t}$ for the case $|\alpha|=n$ even.} The case  $n$ is even is similar. This completes the first induction and yields the proposition.\end{proof}

\begin{corollary} If $\alpha$ is a ribbon and $D$ is a skew diagram then ${\mathfrak s}_{\alpha \bullet  D }={\mathfrak s}_{\alpha \bullet  D ^\circ}.$ \label{dicorrotation} \end{corollary}

\begin{proof} Both cases $|\alpha|$ odd and $|\alpha|$ even follow from Proposition~\ref{diproprottrans}, \eqref{eq:Qrot} and \eqref{dirotation}.\end{proof}

 \comment{For $|\alpha|$ odd,
$${\mathfrak s}_{\alpha\bullet D }={\mathfrak s}_{\alpha^\circ\bullet D }={\mathfrak s}_{(\alpha^\circ\bullet D )^\circ}={\mathfrak s}_{\alpha\bullet D ^\circ}$$
where we used Proposition~\ref{diproprottrans} for the first, \eqref{eq:Qrot} for the second and \eqref{dirotation} for the third equality.

For $|\alpha|$ even,
$${\mathfrak s}_{\alpha\bullet D }={\mathfrak s}_{\alpha^\circ\bullet D }={\mathfrak s}_{(\alpha^\circ\bullet D )^\circ}={\mathfrak s}_{\alpha\bullet( D ^t)^\circ}={\mathfrak s}_{\alpha\bullet D ^\circ}$$
where we used Proposition~\ref{diproprottrans} for the first, \eqref{eq:Qrot} for the second, \eqref{dirotation} for the third and Proposition \ref{diproprottrans} for the fourth equality. }

\begin{corollary}If $\alpha$ is a ribbon and $D$ is a skew diagram then ${\mathfrak s}_{\alpha\bullet  D }={\mathfrak s}_{\alpha^t \bullet  D }$. \label{dicortranspose} \end{corollary}

\begin{proof} Both cases $|\alpha|$ odd and $|\alpha|$ even follow from Proposition~\ref{diproprottrans}, \eqref{eq:Qtr} and \eqref{ditransposition}.\end{proof}
\comment{For $|\alpha|$ even,
$$  {\mathfrak s}_{\alpha\bullet  D }={\mathfrak s}_{(\alpha\bullet  D )^t}={\mathfrak s}_{\alpha^t\bullet  D }  $$
where we used \eqref{eq:Qtr} for the first and \eqref{ditransposition} for the second equality.
 
For $|\alpha|$ odd,
$$  {\mathfrak s}_{\alpha\bullet  D }= {\mathfrak s}_{(\alpha\bullet  D )^t}=  {\mathfrak s}_{\alpha^t\bullet  D ^t}= {\mathfrak s}_{\alpha^t\bullet  D } $$
where we used \eqref{eq:Qtr} for the first, \eqref{ditransposition} for the second and Proposition~\ref{diproprottrans} for the third equality. }

We can also derive new ordinary skew Schur $Q$-function equalities from known ones.

\begin{proposition} For skew diagrams $D$ and $E$, if ${\mathfrak s}_ D ={\mathfrak s}_ E $ then ${\mathfrak s}_{ D \odot  D ^t}={\mathfrak s}_{ D \cdot  D ^t}={\mathfrak s}_{ E \cdot E ^t}={\mathfrak s}_{ E  \odot E ^t} $. \label{diunexplained}\end{proposition}

\begin{proof} Note that $ D \odot  D ^t=2\bullet  D $ and $ D  \cdot  D ^t=11\bullet  D $. Since $2^t = 11$, we have by Corollary \ref{dicortranspose} that ${\mathfrak s}_{ D \odot  D ^t}={\mathfrak s}_{ D \cdot  D ^t}$. The result follows from \eqref{eq:Qmult} with $E= D^t$ yielding
\begin{equation}{\mathfrak s}^2_ D =2{\mathfrak s}_{ D \odot  D ^t} \label{diunexplainedalpha}.
\end{equation}\end{proof}

\comment{Now by \eqref{eq:Qmult} we have
$${\mathfrak s}_ D {\mathfrak s}_{ D ^t}={\mathfrak s}_{ D \cdot  D ^t}+{\mathfrak s}_{ D \odot  D ^t}.$$
Substituting ${\mathfrak s}_ D ={\mathfrak s}_{ D ^t}$, ${\mathfrak s}_{ D \odot  D ^t}={\mathfrak s}_{ D \cdot  D ^t}$ and simplifying we get
\begin{E}{\mathfrak s}^2_ D =2{\mathfrak s}_{ D \odot  D ^t} \label{diunexplainedalpha}.\end{E}
Similarly 

\begin{E}{\mathfrak s}^2_ E =2{\mathfrak s}_{ E \odot  E ^t} \label{diunexplainedbeta}. \end{E}
Since ${\mathfrak s}_ D ={\mathfrak s}_ E $, \eqref{diunexplainedalpha} and \eqref{diunexplainedbeta} yield ${\mathfrak s}_{ D  \odot  D ^t}={\mathfrak s}_{ E  \odot  E ^t}$. The result now follows.} 

\begin{proposition}\label{prop:power2}
For skew diagrams $ D $ and $ E $, ${\mathfrak s}_ D ={\mathfrak s}_ E $ if and only if $${\mathfrak s}_{\underbrace{2\bullet\cdots\bullet 2}_{n}\bullet D }={\mathfrak s}_{\underbrace{2\bullet\cdots\bullet 2}_{n}\bullet E }.$$
\end{proposition}

\begin{proof} This follows from a straightforward application of \eqref{diunexplainedalpha}.\end{proof}
\comment{From \eqref{diunexplainedalpha}, we have
$${\mathfrak s}^2_ D =2{\mathfrak s}_{2\bullet D }$$
therefore 
$${\mathfrak s}^2_{\underbrace{2\bullet\cdots\bullet 2}_{n-1}\bullet D }=2{\mathfrak s}_{2\bullet(\underbrace{2\bullet\cdots\bullet 2}_{n-1}\bullet D )}=2{\mathfrak s}_{\underbrace{2\bullet\cdots\bullet 2}_{n}\bullet D } $$and 
similarly 
$${\mathfrak s}^2_{\underbrace{2\bullet\cdots\bullet 2}_{n-1}\bullet E }=2{\mathfrak s}_{2\bullet(\underbrace{2\bullet\cdots\bullet 2}_{n-1}\bullet E )}=2{\mathfrak s}_{\underbrace{2\bullet\cdots\bullet 2}_{n}\bullet E }. $$ 

Hence, it is straightforward to see
$${\mathfrak s}_{\underbrace{2\bullet\cdots\bullet 2}_{n}\bullet D }={\mathfrak s}_{\underbrace{2\bullet\cdots\bullet 2}_{n}\bullet E } \Longleftrightarrow {\mathfrak s}_{\underbrace{2\bullet\cdots\bullet 2}_{n-1}\bullet D }={\mathfrak s}_{\underbrace{2\bullet\cdots\bullet 2}_{n-1}\bullet E }.$$
Applying the above equivalence repeatedly, we have 
$${\mathfrak s}_{\underbrace{2\bullet\cdots\bullet 2}_{n}\bullet D }={\mathfrak s}_{\underbrace{2\bullet\cdots\bullet 2}_{n}\bullet E } \Longleftrightarrow {\mathfrak s}_ D ={\mathfrak s}_ E .  $$} 

Before we prove our main result on equality we require the following map, which is analogous to the map $\circ s_D$ in \cite[Corollary 7.4]{HDL2}.

\begin{proposition}\label{prop:wdmap} For a fixed skew diagram $D$, the map
\begin{eqnarray*}
 \Q[z_{\alpha}] &\stackrel{(-)\bullet{\mathfrak s}_ D }{\longrightarrow} &\Omega\\
 z_\alpha&\mapsto &{\mathfrak s}_{\alpha\bullet  D } \\
 0&\mapsto &0
\end{eqnarray*}  descends to a well-defined map $\Omega \rightarrow \Omega$. Hence it is well-defined to set
$$\Qr _\alpha \bullet \Qs _D = \Qs _{\alpha \bullet D}$$where we abuse notation by using $\bullet$ for both the map and the composition of transpositions.
\end{proposition}

\begin{proof}
Observe that by Theorem~\ref{ribbonrelations2}  it suffices to prove that the expressions 
 $$z_{\alpha}\ z_{\beta}-z_{\alpha\cdot \beta} - z_{\alpha \odot \beta}$$for ribbons $\alpha, \beta$ and
$$z_{2m} - z_{\underbrace{1\ldots1}_{2m}}$$
for all positive integers $m$, are mapped to 0 by $(-)\bullet{\mathfrak s}_ D$.

For the first expression, observe that for ribbons $\alpha, \beta$ and skew diagram $D$
$$\begin{array}{c}
(\alpha\cdot\beta)\bullet  D =(\alpha\bullet D )\cdot(\beta\bullet D ') \\
(\alpha\odot\beta)\bullet  D =(\alpha\bullet D )\odot(\beta\bullet D ')
\end{array}$$
where $ D '= D $ when $|\alpha|$ is even and $ D '= D ^t$ otherwise. Therefore 
$$z_\alpha z_\beta-z_{\alpha\cdot\beta}-z_{\alpha\odot\beta}$$is mapped to 
$$\begin{array}{ll}
 & {\mathfrak s}_{\alpha\bullet  D }{\mathfrak s}_{\beta\bullet  D }-{\mathfrak s}_{(\alpha\cdot\beta)\bullet D }-{\mathfrak s}_{(\alpha\odot\beta)\bullet D }\\
 = & {\mathfrak s}_{\alpha\bullet  D }{\mathfrak s}_{\beta\bullet  D '}-{\mathfrak s}_{(\alpha\bullet D )\cdot(\beta\bullet D ')}-{\mathfrak s}_{(\alpha\bullet D )\odot(\beta\bullet D ')}\\
 = & 0
 \end{array}$$
where we used the above observation and Proposition~\ref{diproprottrans}  for the first, and \eqref{eq:Qmult} for the second equality.

For the second expression, observe
$$z_{2m}-z_{\underbrace{1\ldots1}_{2m}}$$
goes to 
 $${\mathfrak s}_{2m\bullet  D }-{\mathfrak s}_{\underbrace{1\ldots1}_{2m}\bullet D }= {\mathfrak s}_{2m\bullet  D }-{\mathfrak s}_{2m\bullet  D }=0$$where we used Corollary~\ref{dicortranspose} for the first equality. 
\end{proof}

%We also require

\begin{proposition} For ribbons $\alpha$, $\beta$ and skew diagram $D$, if ${\mathfrak r}_\alpha={\mathfrak r}_\beta$ then ${\mathfrak s}_{\alpha\bullet D }={\mathfrak s}_{\beta\bullet D }$.
\label{dibulletonright}
\end{proposition}

\begin{proof}
This follows by Proposition~\ref{prop:wdmap}.
\comment{ we have
$${\mathfrak s}_{\alpha\bullet D }={\mathfrak r}_\alpha \bullet {\mathfrak s}_ D ={\mathfrak r}_\beta \bullet {\mathfrak s}_ D ={\mathfrak s}_{\beta\bullet D }$$
where the second equality holds because we are given ${\mathfrak r}_\alpha={\mathfrak r}_\beta$. This completes the proof.}
\end{proof}

We now come to our main result on equality of ordinary skew Schur $Q$-functions.

\begin{theorem}\label{the:bigone}
For ribbons $\alpha _1, \ldots , \alpha _m$ and skew diagram $D$ the ordinary skew Schur $Q$-function indexed by
$$\alpha _1 \bullet \cdots \bullet \alpha _m \bullet D$$is equal to the ordinary skew Schur $Q$-function indexed by
$$\beta _1 \bullet \cdots \bullet \beta _m \bullet E$$where
$$ 
\beta _i \in \{ \alpha _i, \alpha _i ^t, \alpha _i ^\circ , (\alpha _i ^t)^\circ = (\alpha _i ^\circ)^t\} \quad 1\leq i \leq m,$$
$$ 
E\in \{ D, D^t, D^\circ , (D^t)^\circ = (D^\circ)^t\}.$$
\end{theorem}

\begin{proof}
We begin by restricting our attention to ribbons and proving that for ribbons $\alpha _1, \ldots , \alpha _m$
$$\Qr _{\alpha _1 \bullet \cdots \bullet \alpha _m} = \Qr _{\beta _1 \bullet \cdots \bullet \beta _m}$$where
$ 
\beta _i \in \{ \alpha _i, \alpha _i ^t, \alpha _i ^\circ , (\alpha _i ^t)^\circ = (\alpha _i ^\circ)^t\} \quad 1\leq i \leq m$.

To simplify notation let $\lambda= \alpha_1 \bullet \cdots \bullet \alpha_m$ and $\mu=\beta_1 \bullet \cdots \bullet \beta_m$ where $\beta_i=\{\alpha_i, {\alpha_i}^t, {\alpha_i}^\circ, ({\alpha_i}^t)^\circ\}$ for $1 \leq i \leq m$.

Let $i$ be the smallest index in $\mu$ such that $\alpha_i \neq \beta_i$. Suppose $\beta_i={\alpha_i}^t$, then by the associativity of $\bullet$
\begin{E} {\mathfrak r}_\mu={\mathfrak r}_{(\alpha_1 \bullet \cdots \bullet \alpha_{i-1}) \bullet({\alpha_i}^t\bullet \beta_{i+1}\bullet \cdots \bullet \beta_m)}={\mathfrak r}_{(\alpha_1 \bullet \cdots \bullet \alpha_{i-1}) \bullet({\alpha_i}^t\bullet \beta_{i+1}\bullet \cdots \bullet \beta_m)^t}={\mathfrak r}_{\alpha_1 \bullet \cdots \bullet \alpha_{i-1} \bullet{\alpha_i}\bullet (\beta_{i+1}\bullet \cdots \bullet \beta_m)'} \label{multitranspose} 
\end{E}where we used Proposition~\ref{diproprottrans} for the second and \eqref{ditransposition} for the third equality. Note that $(\beta_{i+1}\bullet \cdots \bullet \beta_m)'=\beta_{i+1}\bullet \cdots \bullet \beta_m$ if $|\alpha_i|$ is even, and $(\beta_{i+1}\bullet \cdots \bullet \beta_m)'=(\beta_{i+1}\bullet \cdots \bullet \beta_m)^t$ if $|\alpha_i|$ is odd. 

Now suppose $\beta_i={\alpha_i}^\circ$, then
\begin{E} {\mathfrak r}_\mu={\mathfrak r}_{(\alpha_1 \bullet \cdots \bullet \alpha_{i-1}) \bullet({\alpha_i}^\circ \bullet \beta_{i+1}\bullet \cdots \bullet \beta_m)}={\mathfrak r}_{(\alpha_1 \bullet \cdots \bullet \alpha_{i-1}) \bullet({\alpha_i}^\circ\bullet \beta_{i+1}\bullet \cdots \bullet \beta_m)^\circ}={\mathfrak r}_{\alpha_1 \bullet \cdots \bullet \alpha_{i-1} \bullet{\alpha_i}\bullet (\beta_{i+1}\bullet \cdots \bullet \beta_m)'} \label{multirotation} 
\end{E}where we used Corollary~\ref{dicorrotation} for the second and \eqref{dirotation} for the third equality. Note that $(\beta_{i+1}\bullet \cdots \bullet \beta_m)'=(\beta_{i+1}\bullet \cdots \bullet \beta_m)^\circ$ if $|\alpha_i|$ is odd, and $(\beta_{i+1}\bullet \cdots \bullet \beta_m)'=((\beta_{i+1}\bullet \cdots \bullet \beta_m)^t)^\circ$ if $|\alpha_i|$ is even.

For the case $\beta_i=({\alpha_i}^t)^\circ$ we can combine \eqref{multitranspose} and \eqref{multirotation} to arrive at
$${\mathfrak r}_\mu ={\mathfrak r}_{\alpha_1 \bullet \cdots \bullet \alpha_{i-1} \bullet{\alpha_i}\bullet (\beta_{i+1}\bullet \cdots \bullet \beta_m)'} $$and 
$$(\beta_{i+1}\bullet \cdots \bullet \beta_m)'\in \{ (\beta_{i+1}\bullet \cdots \bullet \beta_m), (\beta_{i+1}\bullet \cdots \bullet \beta_m)^t, (\beta_{i+1}\bullet \cdots \bullet \beta_m)^\circ, ((\beta_{i+1}\bullet \cdots \bullet \beta_m)^t)^\circ\}.$$

Iterating the above process for each of the three cases, we recover ${\mathfrak r}_\lambda$. 

Applying Proposition~\ref{dibulletonright}, we have 
$${\mathfrak s}_{\alpha_1\bullet\cdots\bullet\alpha_m\bullet D}={\mathfrak s}_{\beta_1\bullet\cdots\bullet\beta_m\bullet D}.$$
Using Corollary \ref{dicorrotation} and Proposition \ref{diproprottrans} we know that ${\mathfrak s}_{\beta_1\bullet\cdots\bullet\beta_m\bullet D}={\mathfrak s}_{\beta_1\bullet\cdots\bullet\beta_m\bullet E}$ where $E=\{D, D^t,D^\circ,(D^t)^\circ\}$. The assertion follows from combining the latter equality with the above equality.
\end{proof}

%%%%%%%% AT THE END Ribbon Section %%%%%%%%%
\section{\texorpdfstring{Ribbon Schur $Q$-functions}{Ribbon Schur Q-functions}}\label{sec:ribbons}

We have seen that ribbon Schur $Q$-functions yield a natural basis for $\Omega$ in Corollary~\ref{cor:rbasis} and establish a generating set of relations for $\Omega$ in Theorems~\ref{ribbonrelations} and \ref{ribbonrelations2}. Now we will see how they relate to enumeration in graded posets.

Let $NC=\mathbb{Q} \langle y_1, y_2, \ldots \rangle$ be the free associative algebra on countably many generators $y_1, y_2, \ldots$ then \cite{BilleraLiu} showed that $NC$ is isomorphic to the non-commutative   algebra of flag-enumeration functionals on graded posets. Furthermore, they showed that the non-commutative algebra of flag-enumeration functionals on Eulerian posets is isomorphic to
$$A_\mathcal{E} = NC / \langle \chi _2, \chi _4 , \ldots\rangle$$where $\chi _{2m}$ is the even Euler form $\chi _{2m}=\sum _{r+s = 2m} (-1)^{r} y_ry_s$. Given a composition $\alpha = \alpha _1 \alpha _2\cdots \alpha _{\ell(\alpha)}$, the \emph{flag-$f$ operator} $y_\alpha$ is
$$y_\alpha = y _{\alpha _1}y_{\alpha _2} \cdots y_{\alpha _{\ell(\alpha)}}$$and the \emph{flag-$h$ operator}  $\mathfrak{h} _\alpha$ is
$$\mathfrak{h} _\alpha= (-1)^{\ell(\alpha)} \sum _{\beta \coarser \alpha} (-1) ^{\ell(\beta)} y _{\beta}$$and $y _\alpha$ and $\mathfrak{h} _\alpha$ are described as being of Eulerian posets if we view them as elements of $A_\mathcal{E}$.

We can now give the relationship between $A_\mathcal{E}$ and $\Omega$.

\begin{theorem}\label{the:commconnection}
Let $\alpha$ be a composition. The non-commutative analogue of $q_\alpha $ is the flag-$f$operator of Eulerian posets, $y_\alpha$. Furthermore, the non-commutative analogue of ${\mathfrak r} _\alpha$ is the flag-$h$ operator of Eulerian posets, $\mathfrak{h}_\alpha$.
\end{theorem}

\begin{proof}
Consider the map  
\begin{eqnarray*}
\psi: A _{\mathcal{E}}&\rightarrow&\Omega\\
%y_\alpha = y_{\alpha _1}\cdots y_{\alpha _k}&\mapsto & q_\alpha = q_{\alpha _1}\cdots q_{\alpha _k}
y_i&\mapsto &q_i
\end{eqnarray*}
extended multiplicatively and by linearity.

By \cite[Proposition 3.2]{BilleraLiu}  all relations in $A _{\mathcal{E}}$ are generated by all $\chi _{n}= \sum _{i+j=n} (-1)^iy_iy_j$.
Hence $\psi (\chi _n)= \sum _{i+j=n} (-1)^iq_iq_j = 0$ by \eqref{eq:qrels}, and hence $\psi$ is well-defined. 
%Additionally, $\psi (y_\alpha y_\beta)= q_\alpha q_\beta = \psi (y_\alpha)\psi(y_\beta)$ 
Therefore, we have that $\psi$ is an algebra homomorphism.
Since the flag-$h$ operator  of Eulerian posets, $\mathfrak{h}_\alpha$, is defined to be
$$\mathfrak{h}_\alpha= \sum _{\beta \coarser \alpha} (-1)^{\ell(\alpha)-\ell (\beta)} y _\beta $$ we have
$\psi (\mathfrak{h}_\alpha) = {\mathfrak r} _\alpha  $ by \eqref{eq:Qrib}.\end{proof}

\comment{\begin{eqnarray*}\psi (\mathfrak{h}_\alpha) &=& \psi ( \sum _{\beta \coarser \alpha} (-1)^{\ell(\alpha)-\ell (\beta)} y _\beta)\\
&=& \sum _{\beta \coarser \alpha} (-1)^{\ell(\alpha)-\ell (\beta)} q _\beta\\
&=& (-1)^{\ell(\alpha)}\sum _{\beta \coarser \alpha} (-1)^{\ell (\beta)} q _{\lambda(\beta)}\\
&=& {\mathfrak r} _\alpha  \mbox{ by \eqref{eq:Qrib}. } \end{eqnarray*}}

\begin{remark} Note that we have the following commutative diagram
$$\xymatrix{
NC\ar@{->}[r] ^{\theta ^N} \ar@{->>}[d]  _\phi& A_{\mathcal{E}} \ar@{->>}[d] ^\psi \\
\Lambda \ar@{->}[r]^{\theta } & \Omega}$$
where $\phi(y_i)=h_i$ and $h_i$ is the $i$-th homogeneous symmetric function, and $\theta ^N (y_i)=y_i$ is the non-commutative analogue of the   map $\theta$. Abusing notation, and denoting $\theta ^N$ by $\theta$ we summarize the relationships between non-symmetric, symmetric and quasisymmetric functions as follows

$$
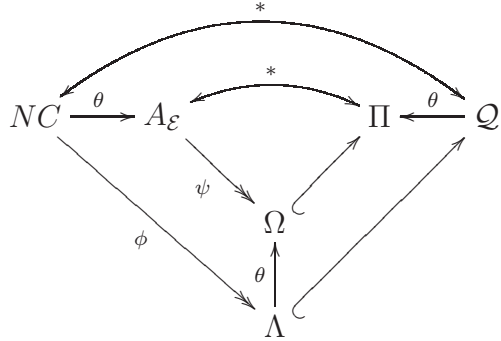
\xymatrix{
NC\ar@{->}[r] ^{\theta } \ar@{->>}[rrdd]  _\phi \ar@/^3pc/@{<->}[rrrr]^\ast & A_{\mathcal{E}} \ar@{->>}[rd] _\psi \ar@/^1pc/@{<->}[rr]^\ast && \Pi & \mathcal{Q}\ar@{->}[l] _\theta\\
&& \Omega \ar@{_{(}->}[ru]\\
&&\Lambda \ar@{->}[u]^{\theta } \ar@{_{(}->}[rruu]
}$$where $\mathcal{Q}$ is the algebra of qusaisymmetric functions and $\Pi$ is the algebra of peak quasisymmetric functions.

For the interested reader, the duality between $NC$ and $\mathcal{Q}$ was established through \cite{Gelfand, Gessel, MR}, and between $A_\mathcal{E}$ and $\Pi$ in \cite{BMSvW}. The commutative diagram connecting $\Omega, \Lambda, \Pi$ and $\mathcal{Q}$ can be found in \cite{Stem}, and the relationship between $NC$ and $\Lambda$ in \cite{Gelfand}.
\end{remark}

\subsection{\texorpdfstring{Equality of ribbon Schur $Q$-functions}{Equality of ribbon Schur Q-functions}} From the above uses and connections it seems worthwhile to restrict our attention to ribbon Schur $Q$-functions in the hope that they will yield some insight into the general solution of when two skew Schur $Q$-functions are equal, as was the case with ribbon Schur functions \cite{HDL, HDL3, HDL2}. Certainly our search space is greatly reduced due to the following proposition.
\begin{proposition}
Equality of skew Schur $Q$-functions restricts to ribbons. That is, if ${\mathfrak r} _\alpha = Q_D$ for a skew diagram $D$ then the shifted skew diagram $\tilde{D}$ must be a ribbon.
\end{proposition}
\begin{proof} Recall  that by definition
$$Q_{D}=\sum _T x^T$$where  the sum is over all weakly amenable tableaux of shape $\tilde{D}$.

If $D$ has $n$ cells,  we now consider the coefficient of $x_1 ^n$ in three scenarios.
\begin{enumerate}
\item $\tilde{D}$ is a ribbon: $[Q_D] _{x_1^n}=2$, which arises from the weakly amenable tableaux where every cell that has a cell to its left must be occupied by
$1$, every cell that has a cell below it must be occupied by $1'$, and the
bottommost and leftmost cell can be occupied by either $1$ or $1'$.
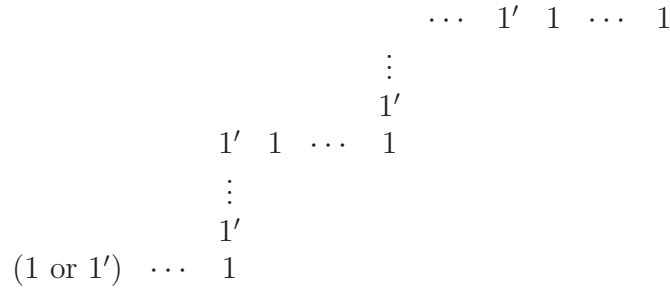
$$\begin{matrix}
&&&&&&\cdots&1'&1&\cdots&1\\
&&&&&\vdots\\
&&&&&1'\\
&&1'&1&\cdots&1\\
&&\vdots\\
&&1'\\
(1\mbox{ or }1')&\cdots & 1\\
\end{matrix}$$
\item $\tilde{D}$ is disconnected and each connected component is a ribbon: $[Q_D] _{x_1^n}=2^{c}$ where $c$ is the number of connected components. This is because the leftmost cell in the bottom row of all  components can be filled with $1$ or $1'$ to create a weakly amenable tableau, and the remaining cells of each connected component can be filled as in the last case.
\item $\tilde{D}$ contains a $2\times 2$ subdiagram: $[Q_D] _{x_1^n}=0$ as the $2\times 2$ subdiagram cannot be filled only with $1$ or $1'$ to create a weakly amenable tableau.
\end{enumerate}

Now note that if ${\mathfrak r} _\alpha = Q_D$ then the coefficient of $x_1 ^n$ must be the same in both
${\mathfrak r} _\alpha$ and  $Q_D$. From the above case analysis we see that the coefficient of $x_1 ^n$ in ${\mathfrak r} _\alpha$
is 2, and hence also in $Q_D$. Therefore, by the above case analysis, $\tilde{D}$ must also be a ribbon.
\end{proof}

We now recast our main results from the previous section in terms of ribbon Schur $Q$-functions, and use this special case to illustrate our results.

\begin{proposition}
For ribbons $\alpha$ and $\beta$, ${\mathfrak r}_\alpha ={\mathfrak r}_\beta $ if and only if $${\mathfrak r}_{\underbrace{2\bullet\cdots\bullet 2}_{n}\bullet \alpha }={\mathfrak r}_{\underbrace{2\bullet\cdots\bullet 2}_{n}\bullet \beta }.$$ \label{prop:twos}
\end{proposition}

\begin{example} If we know $\Qr _{2\bullet 2 \bullet 2} = \Qr _{3311}= \Qr _{1511} = \Qr _{2\bullet 2 \bullet 11}$ then we have $\Qr _2 = \Qr _{11}$. This would be an alternative  to deducing this result from \eqref{eq:Qtr}.

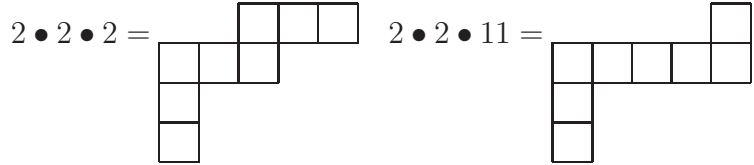
$$2\bullet 2\bullet 2 = \tableau{&&\cell&\cell&\cell\\ \cell&\cell&\cell\\\cell \\ \cell} \quad 2\bullet 2\bullet 11 = \tableau{&&&&\cell \\\cell&\cell&\cell& \cell&\cell\\ \cell \\ \cell}
$$
\end{example}

\begin{remark} Note that the factor 2 appearing in the above proposition is of some fundamental importance since $\Qr _{21\circ 14} = \Qr _{12\circ 14}$ but $\Qr _{3\bullet (21\circ 14)} \neq \Qr _{3\bullet(12\circ 14)}$.
\end{remark}

\begin{proposition}
For ribbons $\alpha, \beta, \gamma$, if ${\mathfrak r}_\alpha={\mathfrak r}_\beta$ then ${\mathfrak r}_{\alpha\bullet \gamma }={\mathfrak r}_{\beta\bullet \gamma }$.\label{prop:ribtorib}
\end{proposition}

\begin{example} Since $\Qr _3 = \Qr _{111}$ by \eqref{eq:Qtr} we have $\Qr _{33141} = \Qr _{3\bullet 31} = \Qr _{111\bullet 31} = \Qr _{3121131}.$

$$3\bullet 31 = \tableau{&&&&&\cell&\cell&\cell\\
&&&\cell&\cell&\cell\\
&&&\cell\\
\cell&\cell&\cell&\cell\\\cell}\quad 111\bullet 31 = \tableau{&&&\cell&\cell&\cell\\&&&\cell\\&&\cell&\cell\\&&\cell\\&&\cell\\\cell&\cell&\cell\\\cell}$$
However, we could also have deduced $\Qr _{33141}  = \Qr _{3121131}$ from the following theorem.

\end{example}

\begin{theorem}
For ribbons $\alpha _1, \ldots , \alpha _m$ the ribbon Schur $Q$-function indexed by
$$\alpha _1 \bullet \cdots \bullet \alpha _m$$
is equal to the 
ribbon Schur $Q$-function indexed by
$$\beta _1 \bullet \cdots \bullet \beta _m$$where
$$\beta _i \in \{ \alpha _i, \alpha _i ^t, \alpha _i ^\circ , (\alpha _i ^t)^\circ = (\alpha _i ^\circ)^t\} \quad 1\leq i \leq m.$$
\end{theorem}

\begin{example}
If $\alpha _1=2$ and $\alpha _2 = 21$ then
$$\Qr _{231} = \Qr _{2121} = \Qr _{132} = \Qr _{1212}$$as
$$2\bullet 21 = \tableau{&&\cell&\cell\\\cell&\cell&\cell\\\cell}\ , 
2^t\bullet 21 = \tableau{&\cell&\cell\\&\cell\\\cell&\cell\\\cell}\ ,
2\bullet (21)^\circ = \tableau{&&&\cell\\&\cell&\cell&\cell\\\cell&\cell}\ ,
2^t\bullet (21)^\circ = \tableau{&&\cell\\&\cell&\cell\\&\cell\\\cell&\cell}\ ,$$but we could have equally well just chosen
$\alpha = 231$ and concluded again
\begin{eqnarray*}
\Qr _{231} &=& \Qr _{(231)^t} = \Qr _{(231)^\circ} = \Qr _{((231)^t )^\circ}\\
&=&\Qr _{2121} = \Qr _{132} = \Qr _{1212}.\\
\end{eqnarray*}
\end{example}

We begin to draw our study of ribbon Schur $Q$-functions to a close with the following conjecture, which we prove in one direction, and has been confirmed for ribbons with up to 13 cells.

\begin{conjecture}
For ribbons $\alpha, \beta$ we have $\Qr  _\alpha = \Qr _\beta$ if and only if there exists $j, k, l$  so that
$$\alpha = \alpha _1 \bullet \cdots \bullet \alpha _j \bullet (\gamma _1 \circ \cdots \circ \gamma _k)\bullet \varepsilon _1 \bullet \cdots \bullet\varepsilon _\ell$$and
$$\beta = \beta _1 \bullet \cdots \bullet \beta _j \bullet (\delta _1 \circ \cdots \circ \delta _k)\bullet \eta _1 \bullet \cdots  \bullet \eta _\ell$$where
$$\alpha _i, \beta _i \in \{ 2, 11\}\quad 1\leq i \leq j,$$$$ \delta _i \in \{\gamma _i , \gamma _i ^\circ\} \quad 1\leq i \leq k,$$$$ 
\eta _i \in \{ \varepsilon _i, \varepsilon _i ^t, \varepsilon _i ^\circ , (\varepsilon _i ^t)^\circ = (\varepsilon _i ^\circ)^t\} \quad 1\leq i \leq \ell.$$
\end{conjecture}

To prove one direction note that certainly if $\alpha$ and $\beta$ satisfy the criteria then $\Qr _\alpha = \Qr _\beta$ since by applying $\theta$ to \cite[Theorem 4.1]{HDL} we have
$$\Qr _{\gamma _1 \circ \cdots \circ \gamma _k} = \Qr _{\delta _1 \circ \cdots \circ \delta _k}.$$By Proposition~\ref{prop:twos} and Corollary~\ref{dicortranspose} we get
$$\Qr _{11\bullet (\gamma _1 \circ \cdots \circ \gamma _k)} = \Qr _{2\bullet (\gamma _1 \circ \cdots \circ \gamma _k)}= \Qr _{2\bullet (\delta _1 \circ \cdots \circ \delta _k)}= \Qr _{11\bullet (\delta _1 \circ \cdots \circ \delta _k)}$$and performing this repeatedly we get 
$$\Qr _{\alpha _1 \bullet \cdots \bullet \alpha _j \bullet (\gamma _1 \circ \cdots \circ \gamma _k)}= \Qr _{\beta _1 \bullet \cdots \bullet \beta _j\bullet (\delta _1 \circ \cdots \circ \delta _k)}.$$By Proposition~\ref{prop:ribtorib}, Proposition~\ref{diproprottrans} and Corollary~\ref{dicorrotation} we get 
\begin{eqnarray*}\Qr _{\beta _1 \bullet \cdots \bullet \beta _j\bullet (\delta _1 \circ \cdots \circ \delta _k)\bullet\varepsilon _1} &=&
\Qr _{\alpha _1 \bullet \cdots \bullet \alpha _j \bullet (\gamma _1 \circ \cdots \circ \gamma _k)\bullet\varepsilon _1}\\
&=&
\Qr _{\alpha _1 \bullet \cdots \bullet \alpha _j \bullet (\gamma _1 \circ \cdots \circ \gamma _k)\bullet\varepsilon _1^t}\\
&=&
\Qr _{\alpha _1 \bullet \cdots \bullet \alpha _j \bullet (\gamma _1 \circ \cdots \circ \gamma _k)\bullet\varepsilon _1^\circ}\\
&=&
\Qr _{\alpha _1 \bullet \cdots \bullet \alpha _j \bullet (\gamma _1 \circ \cdots \circ \gamma _k)\bullet(\varepsilon _1^t)^\circ}
\end{eqnarray*}and performing this repeatedly and noting the associativity of $\bullet$ we obtain one direction of our conjecture.

Proving the other direction may be difficult, as a useful tool in studying equality of skew Schur functions was the irreducibility of those indexed by a connected skew diagram \cite{HDL2}. However, irreducibility is a more complex issue when studying the equality of skew Schur $Q$-functions, as illustrated by restricting to ribbon Schur $Q$-functions.

\begin{proposition}\label{prop:irrrib} Let $\alpha$ be a ribbon
\begin{enumerate}
\item for $|\alpha|$ odd,  ${\mathfrak r}_\alpha$ is irreducible 
\item for $|\alpha|$ even, there are infinitely many examples in which ${\mathfrak r}_\alpha$ is irreducible and infinitely many examples in which ${\mathfrak r}_\alpha$ is reducible
\end{enumerate}
considered as an element of ${\mathbb Z}[q_1,q_3,\ldots]$. 
\end{proposition}

\begin{proof} We first prove the first assertion. Let $|\alpha|=n$, where $n$ is an odd integer. Using \eqref{eq:Qrib}, we have 
$${\mathfrak r}_\alpha= \pm q_n+r$$
in which $r$ involves only $q_1,q_3,\ldots,q_{n-2}$. This shows that ${\mathfrak r}_\alpha$ is irreducible in ${\mathbb Z}[q_1,q_3,\ldots]$.

For the second assertion, note that 
$$\Qr _\alpha ^2 = 2\Qr _{\alpha \odot \alpha ^t}$$by \eqref{diunexplainedalpha}.  Hence, ${\mathfrak r}_{\alpha\odot\alpha^t}$ is reducible for every choice of $\alpha$. Further, we show that ${\mathfrak r}_{(4x)2}$ is irreducible in ${\mathbb Z}[q_1,q_3,\ldots]$ for every positive integer $x$. By \eqref{eq:Qrib} we have 
$${\mathfrak r}_{(4x)2}=q_{4x}q_2-q_{4x+2}=-q_{4x+1}q_1\underbrace{+2q_{4x}q_2-q_{4x-1}q_3+\cdots+q_{2x+2}q_{2x}}_{A}-\frac{q_{2x+1}^2}{2} $$
where we substituted $q_{4x+2}$ using \eqref{qrelations} and simplified for the second equality. We use \eqref{qrelations} to reduce the terms in part $A$ into $q$'s with odd subscripts; however, note that no term in $A$ would contain $q_{4x+1}$ and the terms that contain $q_{2x+1}$ have at least two other $q$'s in them. Since the expansion of ${\mathfrak r}_{(4x)2}$ has $-q_{4x+1}q_1$ and no other term containing $q_{4x+1}$, if ${\mathfrak r}_{(4x)2}$ is reducible then $q_1$ has to be a factor of it. But because we have a non-vanishing term $\frac{q_{2x+1}^2}{2}$ in the expansion of ${\mathfrak r}_\alpha$, $q_1$ cannot be a factor.\end{proof}

\end{document}